\documentclass[12pt,a4paper]{article}
\usepackage[T2A]{fontenc}
\usepackage[english]{babel}
\usepackage{srcltx}\usepackage{graphicx,color}
\usepackage[all]{xy}
\usepackage{latexsym,amsfonts,amssymb,amsmath,longtable,
mathrsfs,theorem,cite
}
\usepackage{stmaryrd} 
\usepackage{calrsfs}
\usepackage
{hyperref}
\renewcommand{\leq}{\leqslant}
\renewcommand{\geq}{\geqslant}

\newcommand{\Oo}{\mathrm O}
\newcommand{\F}{\mathbb F}

\newcommand{\D}{\mathscr D}
\newcommand{\E}{\mathscr E}
\newcommand{\X}{\mathfrak X}
\newcommand{\G}{\mathfrak G}\newcommand{\qed}{\hfill{$\Box$}}
\theoremstyle{plain}
{

\newtheorem{Lemma}{{\bfseries Lemma}}[section]
{\theorembodyfont{\itshape}

\newtheorem{Theo}{{\bfseries Theorem}}

\newtheorem{Probl}{{\bfseries Problem}}
\newtheorem{Cor}{{\bfseries Corollary}}}
\theorembodyfont{\rm}

\newtheorem{Def}{{\bfseries Definition}}

\DeclareMathOperator{\Inn}{Inn} 
 \DeclareMathOperator{\GL}{GL}
\DeclareMathOperator{\Aut}{Aut} \DeclareMathOperator{\Out}{Out} \DeclareMathOperator{\SL}{SL}\DeclareMathOperator{\PSL}{PSL}
\DeclareMathOperator{\Sym}{Sym} 
\DeclareMathOperator{\Sp}{Sp} \DeclareMathOperator{\GU}{GU}
\DeclareMathOperator{\Alt}{Alt}
 \DeclareMathOperator{\Hall}{Hall}
\DeclareMathOperator{\PSp}{PSp}
\renewcommand{\k}{\mathrm{k}} \DeclareMathOperator{\sm}{sm} \DeclareMathOperator{\m}{m}

\newcommand{\PGL}{{\mathrm{PGL}}}
\newcommand{\PGU}{{\mathrm{PGU}}}
\newcommand{\SU}{{\mathrm{SU}}}
\renewcommand{\simeq}{\cong}
\newcommand{\PSU}{{\mathrm{PSU}}}

\renewcommand{\P}{\mathrm{P}}
\newcommand{\splitext}{\,\colon\!}
\newcommand{\arbitraryext}{\,\ldotp}
\newcommand{\nonsplitext}{\,{}^{\text{\normalsize{\textperiodcentered}}}}

\renewcommand{\le}{\leqslant}
\renewcommand{\ge}{\geqslant}
\setlength{\textwidth}{160mm}
\setlength{\textheight}{250mm} \headheight14mm \headsep0mm
\oddsidemargin0mm \topmargin -10mm

\title{\vspace{-4cm} \hfill{\normalsize 20D60}{
\fontfamily{cmr} \fontseries{bx} \selectfont \\ \vspace{1cm}
The reduction theorem\\ for relatively maximal subgroups}
\thanks{ Funding: The first  author is supported by
a NNSF grant of China 
and Wu Wen-Tsun Key Laboratory of Mathematics of Chinese Academy of Sciences.
The second  author is supported by Chinese Academy of Sciences President's International Fellowship Initiative, PIFI, (Grant \#~2016VMA078), and by  by RFBR and BRFBR, project number  20-51-00007. The second and third authors  are supported by  the Program of Fundamental Scientific Research of the SB RAS No. I.1.1., project  number 0314-2016-0001. }}
\date{}

\author{Wenbin Guo\\
{\small 1. School of Science, Hainan University
Haikou, }\\ {\small Hainan, 570228, P.R. China,   and}\\
{\small 2. Department of Mathematics, University of Science and
Technology of China,}\\ {\small Hefei 230026, P. R. China}\\
{\small E-mail:
wbguo@ustc.edu.cn}\\ \\
Danila O. Revin, and Evgeny P. Vdovin\\
{\small 1. Sobolev Institute of Mathematics SB RAS, and}\\
{\small 2.  Novosibirsk State University,}\\
{\small Novosibirsk 630090, Russia}\\
{\small E-mail: revin@math.nsc.ru,  vdovin@math.nsc.ru}
}

\begin{document}

\maketitle
\pagenumbering{arabic}
\begin{abstract}

Let $\X$ be a class of finite groups closed under taking subgroups,
homomorphic images and extensions. It is known that if $A$ is a normal subgroup of a finite group $G$ then the image of an $\X$-maximal subgroup $H$ of $G$ in $G/A$ is not, in general, $\X$-maximal in~$G/A$. We say that the reduction $\X$-theorem holds for a finite group $A$ if, for every finite group $G$ that is an extension of~$A$ (i.\,e. contains $A$ as a normal subgroup), the number of conjugacy classes of $\X$-maximal subgroups in $G$ and $G/A$ is the same. The reduction $\X$-theorem for $A$ implies that $HA/A$ is $\X$-maximal in $G/A$ for every extension $G$ of $A$ and every $\X$-maximal subgroup $H$ of~$G$. In this paper, we prove that the reduction $\X$-theorem holds for $A$ if and only if all $\X$-maximal subgroups are conjugate in $A$ and classify the finite groups with this property in terms of composition factors. 

\smallskip

\end{abstract}

\hfill{\em To the 110th anniversary of the birth of Helmut Wielandt.}

\section{Introduction}

\subsection{ Relatively ma\-xi\-mal subgroups and reduction theorems}
In this paper, we consider only finite groups. $G$ always denotes a finite
group.

Since its inception in the papers by \'E.\,Galois \cite{Galois} and C.\,Jordan \cite{Jordan0,Jordan}, group theory has had the following as one of its central problems: {\it given a group $G$ and a class $\X$ of groups (the class of solvable groups, for example), determine the $\X$-subgroups of $G$} (i.\,e. the subgroups belonging to $\X$). If ${\mathfrak X}$ has good properties resembling those of solvable groups
then to solve the general problem it suffices to know the so-called {\it maximal $\mathfrak{X}$-subgroups} (or {\em $\mathfrak{X}$-maximal subgroups}), i.\,e. the maximal by inclusion $\mathfrak{X}$-subgroups\footnote{To distinguish $\X$-maximal subgroups for a class $\X$ from maximal (among proper) ones, Wielandt~\cite{WieCanb} suggested using the term ``{\it relatively maximal}'' (w.\,r.\,t. a class $\X$). }.
Following H.\,Wielandt \cite{Wie3,Wie4}, a nonempty class $\X$ is said to be {\it complete} ({`vollst\"{a}ndig}') if subgroups and homomorphic images of an $\X$-group and extensions of an $\X$-group by an $\X$-group are always $\X$-groups\footnote{Recall that a group $G$ is an {\it extension}
of a group $A$ by a group $B$ if there is an epimorphism $G\rightarrow B$
with kernel isomorphic to $A$.}.  Examples of complete classes are
\begin{itemize}
  \item $\mathfrak{G}$, the class of all finite groups;
  \item $\mathfrak{S}$, the class of all finite solvable groups.
  \item $\mathfrak{G}_\pi$, the class of all $\pi$-groups for a set $\pi$ of
  primes (i.\,e. the groups $G$ such that every prime divisor of $|G|$ belongs to~$\pi$);
  \item $\mathfrak{S}_\pi$, the class of all solvable $\pi$-groups for a set
  $\pi$ of primes.
\end{itemize}
In fact, the latter two cases are extremal for every complete $\X$: if we denote by
$\pi=\pi(\X)$ the union of the sets $\pi(G)$ of prime divisors of $|G|$, where $G$
runs  over~$\X$, then 
$$
\mathfrak{S}_\pi\subseteq\X\subseteq\mathfrak{G}_\pi.
$$

We will henceforth assume that $\mathfrak{X}$ is a fixed complete class. We denote by $\m_\X(G)$ the set of {\it $\X$-ma\-xi\-mal} subgroups of $G$.

Following an approach practiced in group theory, it is natural to try to reduce the problem of determining   $\X$-ma\-xi\-mal subgroups of $G$ to the sections of a normal or subnormal series $$G=G_0\ge G_1\ge\dots\ge G_n=1,$$ i.\,e. a series such that, for  each $i=1,\dots,n$, $G_i$ is normal in $G$ or $G_{i-1}$, respectively.  In this case, instead of a subgroup $H$ of $G$,  the {\it projection} $$H^i=(H\cap G_{i-1})G_i/G_i$$ in  each section $G^i=G_{i-1}/G_i$ is considered. It is easy  to see that $H\in \X$ if and only if $H^i\in \X$ for all $i=1,\dots,n$. Moreover, if $H^i\in \m_\X(G^i)$ for all $i=1,\dots,n$, then $H\in \m_\X(G)$~\cite[(2.5) and~(2.6)]{WieCanb}.  In order to apply  inductive arguments to studying $\X$-maximal subgroups, it is natural to ask: is every projection $H^i$ an $\X$-maximal subgroup of $G^i$ if $H$ is $\X$-maximal in $G$? This question is equivalent to the following.
\begin{itemize}
\item[$(*)$] Assume that $A$ is a normal subgroup and $H$ is an $\X$-maximal subgroup of $G$. Is it true that $HA/A\in \m_\X(G/A)$ and $H\cap A\in \m_\X(A)$?
\end{itemize}
Complete classes $\X$ satisfying $(*)$ are known. The class $\mathfrak{G}_p$ of  all $p$-groups for a prime $p$ is one of them: the $\mathfrak{G}_p$-maximal subgroup are exactly the Sylow $p$-subgroups by the Sylow theorem, and, for any Sylow $p$-subgroup $H$ of $G$, if $A$ is normal in $G$ then $H\cap A$ is a Sylow $p$-subgroup of $A$, as well as $HA/A$ is a Sylow $p$-subgroup of $G/A$. More generally, if the index of an $\X$-subgroup $H$ of $G$ is not divisible by the primes in $\pi=\pi(\X)$ then $H\in\m_\X(G)$ and, moreover, $HA/A\in \m_\X(G/A)$ and $H\cap A\in \m_\X(A)$ for each normal subgroup $A$ of $G$. In this case, $H$ is a so-called {\it $\pi$-Hall subgroup}, i.\,e. a $\pi$-subgroup whose index is divisible by no primes from~$\pi$. But, in general, the answer to question $(*)$ is negative.  We are concerned with this question for $H\cap A$ in detail in~\ref{sub4}. Here we consider the first part of this question concerning with the behavior of $\X$-maximal subgroups  under epimorphisms.

Wielandt notices \cite[14.2]{Wie4}, \cite[4.3]{Wie3} that, if there exists a group  $L$ with at least two conjugacy classes of
$\X$-ma\-xi\-mal subgroups\footnote{The existence of a group with at least two conjugacy classes of $\X$-ma\-xi\-mal subgroups for a complete~$\X$ is equivalent to the fact that $\X$ differs from each of the classes $\mathfrak{G}$,  $\mathfrak{G}_\varnothing$, and  $\mathfrak{G}_p$, where $p$ is prime, see~\cite{Epimax}.} for a given  complete class $\X$, then there are no restrictions for an $\X$-subgroup of an arbitrary group $G_0$ to coincide with the the image of an $\X$-maximal subgroup under an appropriate epimorphism onto~$G_0$:
\begin{itemize}
    \item {\it every} (not only $\X$-maximal) $\X$-subgroup of $G_0$ is the image of an
$\X$-ma\-xi\-mal subgroup of the regular wreath product $G=L\wr G_0$ under the natural epimorphism $G\rightarrow G_0$.
\end{itemize}

However, Wielandt  \cite[Sections 12 and~15]{Wie4} shows that, under some restrictions on the kernel of a homomorphism $\phi$ from a group $G$,  there is a natural bijection between the conjugacy classes of $\X$-maximal subgroups of $G$ and those of $G^\phi$.

Denote by $\k_\X(G)$ the  number of conjugacy classes of $\X$-maximal subgroups of~$G$.

\begin{Def} We say that {\it the $\X$-Reduktionssatz $($the reduction $\X$-theorem$)$ holds for a group~$A$}, if $\k_\X(G)= \k_\X(G/N)$ for every group $G$ containing a normal subgroup $N$ isomorphic to~$A$.
\end{Def}

It is not hard to show that if ${\phi:G\rightarrow G_0}$ is an epimorphism then, for every ${K\in\m_\X(G_0)}$, there exists ${H\in\m_\X(G)}$ such that $K=H^\phi$. As a consequence,
 \begin{itemize}
   \item $\k_\X(G_0)\le \k_\X(G)$  and
   \item $\k_\X(G_0)= \k_\X(G) \text{ implies }\m_\X(G_0)=\left\{H^\phi\mid H\in\m_\X(G)\right\}.$
 \end{itemize}
It means that {\it the $\X$-Reduktionssatz for $A$ is equivalent to the following statement}:
 \begin{itemize}
\item
\noindent {\it If $G$ contains a normal subgroup $N$ isomorphic to~$A$ then $H\mapsto HN/N$ maps $\m_\X(G)$ onto $\m_\X(G/N)$ and induces the natural bijection between the conjugacy classes of $\X$-maximal subgroups of $G$ and~$G/N$.}
\end{itemize}
If a group $G$ has a normal subgroup $N$ satisfying the $\X$-Reduktionssatz then we can study $\X$-maximal subgroups in $G/N$ instead of those in~$G$. Thus, in order to determine $\X$-maximal subgroups of finite groups up to conjugation, it would be useful to know all groups with the $\X$-Reduktionssatz. The program of studying such groups was outlined by Wielandt in his lectures in T\"{u}bingen in 1963-64 \cite[Chapter~III]{Wie4}. The goal of this paper is to classify all groups with the $\X$-Reduktionssatz for any complete~$\X$. As a matter of fact, this result is a realization of Wielandt's program.

Wielandt shows \cite[12.9]{Wie4} that the $\X$-Reduktionssatz holds for the so-called {\it $\X$-separable} (`$\X$-reihig') groups, i.\,e. for groups possessing a subnormal series whose sections either are $\X$-groups or have no nontrivial $\X$-subgroups. In particular, the $\X$-Reduktionssatz holds for the solvable groups.

Notice that if the $\X$-Reduktionssatz holds for a group $G$ then ${\k_\X(G)= \k_\X(G/G)=1}$, which means that all $\X$-maximal subgroups are conjugate in $G$. Does the converse statement hold, i.\,e. does the  $\X$-Reduktionssatz hold for $G$ if the $\X$-maximal subgroup of $G$ are conjugate? In another form, this question, as well as the problem of classification of groups with the $\X$-Reduktionssatz, was asked by Wielandt in~\cite{Wie4}, see~\ref{sub2}.

The following statement answers this question and completely determines the groups satisfying the $\X$-Re\-duk\-tions\-satz in terms of their composition factors\footnote{Recall that the {\it composition factors} of a finite group are the sections if its {\it composition series}, i.\,e. a subnormal series whose sections are simple. According to the Jordan--H\"{o}lder theorem, the multiset of composition factors of a group $G$ is an invariant of $G$ up to isomorphisms and does not depend on the composition series.}:

\begin{Theo}\label{Reduktionssatz_main} {\rm (mod CFSG)}
For a finite group $G$, the following statements are equivalent:

\begin{itemize}
  \item[$(i)$] the  $\X$-Reduktionssatz hold for $G$;
  \item[$(ii)$] the  $\X$-maximal subgroups of $G$ are conjugate;
  \item[$(iii)$] for every composition factor $S$ of $G$, either $S\in\X$ or the pair $(S,\X)$ satisfies one of Conditions {\rm I--VII} in Appendix.
\end{itemize}
\end{Theo}

Theorem~\ref{Reduktionssatz_main} can be considered as the most complete among other reduction theorems (see~\cite[sections 12 and~15]{Wie4}, for example) and includes them as particular cases.  Our proof is based on Wielandt's ideas and implements his strategy suggested in~\cite{Wie4}. To prove Theorem~\ref{Reduktionssatz_main}, we have to solve some problems posed by Wielandt, see~\ref{sub2}.

Other ingredients of the proof are connected with Hall subgroups and the so-called $\D_\pi$-groups,  i.\,e.  groups such that the complete analog of the Sylow theorem holds for their $\pi$-subgroups. It is easy to show that if the  $\X$-maximal subgroups are conjugate in $G$ then every  $\X$-maximal subgroup is a $\pi$-Hall subgroup for $\pi=\pi(\X)$ and $G$ is like a $\D_\pi$-group. Hall subgroups and $\D_\pi$-groups are well-studied \cite{GRV,MR,MRV,DpiCl,R2,R4,R5,HallChev,Hall3',NumbCl,Thomp,Surv}. Theorem~\ref{Reduktionssatz_main} is  an application of these results, see~\ref{sub3}.

\subsection{Wielandt's strategy and $\X$-submaximal subgroups}\label{sub2}

In order to prove the equivalence of $(i)$ and $(ii)$ in Theorem~\ref{Reduktionssatz_main}, we have to prove Wielandt's conjecture that the conjugacy of the
$\X$-ma\-xi\-mal subgroups  implies the conjugacy of
the so-called $\X$-sub\-ma\-xi\-mal subgroups introduced in \cite[Definition~15.1]{Wie4} and \cite[p.~170]{Wie3}. We discuss the concept of an $\X$-submaximal subgroup generalizing that of an $\X$-maximal subgroup below. Briefly,  Wielandt proved \cite[15.4]{Wie4} that, for every group, the conjugacy of the $\X$-submaximal subgroups implies the $\X$-Reduktionssatz. As we have seen, this implies the conjugacy of the $\X$-maximal subgroups:
$$
\left[\begin{array}{c}
        \text{the conjugacy of} \\
       \text{the } \X\text{-submaximal} \\
        \text{subgroups}
      \end{array}
\right]\Rightarrow\left[\begin{array}{c}
         \\
       \text{the } \X\text{-Reduktionssatz} \\
        ~
      \end{array}\right]\Rightarrow\left[\begin{array}{c}
        \text{the conjugacy of} \\
       \text{the } \X\text{-maximal} \\
        \text{subgroups}
      \end{array}\right]
$$
The problem of the equivalence of $(i)$ and $(ii)$ in Theorem~\ref{Reduktionssatz_main} may be extended to
\begin{Probl}\label{WieProb1}  {\rm (H.~Wielandt,
\cite[offene Frage zu 15.4]{Wie4})  Does the conjugacy of the
$\X$-ma\-xi\-mal subgroups of a group $G$ imply the conjugacy of
the $\X$-sub\-ma\-xi\-mal subgroups of~$G$? }
\end{Probl}

$$
\left[\begin{array}{c}
        \text{the conjugacy of} \\
       \text{the } \X\text{-submaximal} \\
        \text{subgroups}
      \end{array}
\right]
\stackrel{\mbox{\,?} }{ \Leftarrow }
\left[\begin{array}{c}
        \text{the conjugacy of} \\
       \text{the } \X\text{-maximal} \\
        \text{subgroups}
      \end{array}\right]
$$
This problem is central in our paper and we solve it in the positive.

Wielandt's definitions of $\X$-subma\-xi\-mal subgroups given in \cite[Definition~15.1]{Wie4} and \cite[p.~170]{Wie3} are different and even non-equivalent \cite[Section~2]{RSV}. Although we solve the problems presented in the lectures~\cite{Wie4}, the working definition for us is the broader definition of an $\X$-ma\-xi\-mal subgroup given in the talk\footnote{For the role of $\X$-ma\-xi\-mal subgroups in Wielandt's program presented in~\cite{Wie3}, see  \ref{sub4}.}~\cite{Wie3}.
Following~\cite{RSV},  the $\X$-submaximal subgroups in the sense of \cite[Definition~15.1]{Wie4} are called  strongly $\X$-submaximal.

\begin{Def}\label{Submax def} A subgroup $H$ of a group $G$ is called {\it strongly $\X$-sub\-ma\-xi\-mal} or {\it $\X$-sub\-ma\-xi\-mal
in the sense of}~\cite{Wie4} (respectively, {\it $\X$-submaximal} or {\it $\X$-submaximal in the sense of}~~\cite{Wie3}), if there is an embedding $\phi:G\hookrightarrow G^*$ in a  group
$G^*$ such that $G^\phi$ is normal (respectively, subnormal, i.\,e. a member of a subnormal series) in $G^*$ and $H^\phi=K\cap G^\phi$ for an $\X$-ma\-xi\-mal
subgroup $K$ of $G^*$.
\end{Def}

We denote the set of strongly  $\X$-sub\-ma\-xi\-mal subgroups of~$G$
by~$\sm_\X^\circ(G)$ and the set of $\X$-sub\-ma\-xi\-mal subgroups of~$G$
by~$\sm_\X(G)$. It is clear that
$$
\m_\X^{\phantom{\circ}}(G)\subseteq \sm_\X^{{\circ}}(G)\subseteq  \sm_\X^{\phantom{\circ}}(G).
$$
Thus, if we prove that the conjugacy of the members of $\m_\X(G)$ always implies that of $\sm_\X(G)$, then we solve Problem~\ref{WieProb1} to be understood in either sense.

Let us assume we have already done that. Then, applying induction and the $\X$-Re\-duk\-tions\-satz, one can show that the $\X$-maximal (and, consequently, the $\X$-submaximal) subgroups are conjugate in $G$ if they are conjugate in every composition factor of $G$, cf. \cite[p.~57, Satz~1(a)]{Wie_diary}. Using his result \cite{Wie1}, Wielandt notices that, if  $G$ itself \cite[15.6]{Wie4} or even every composition factor of $G$  \cite[15.4$''$]{Wie4} possesses a nilpotent $\pi$-Hall
subgroup for $\pi=\pi(\X)$, then $G$ has exactly one conjugacy class of (strongly)
$\X$-sub\-ma\-xi\-mal subgroups. In view of this fact, he  asks \cite[offene
Frage, p.~37~(643)]{Wie4}, \cite[p.~57, Frage~2]{Wie_diary}: for  a simple\footnote{In Wielandt's lectures, the word ``simple'' is missing from the wording of this question. But it is clear from the context that it is a misprint. This follows from the meaning of the statement \cite[15.4$''$, p.~37~(643)]{Wie4} which immediately precedes the formulation of \cite[offene
Frage, p.~37~(643)]{Wie4}. The problem of the existence of a {\it simple} group possessing exactly one conjugacy class of $\X$-submaximal subgroups, but containing no nilpotent $\pi$-Hall subgroup, is written down in Wielandt's scientific diary on the last page of the main part \cite[p.~57, Frage~2 (entry dated August 20, 1982)]{Wie_diary}. Examples of nonsimple groups without nilpotent $\pi$-Hall subgroups but with a unique conjugacy class of $\X$-sub\-ma\-xi\-mal subgroups were known to Wielandt, see \cite[12.9]{Wie4}.} group $G\notin\X$, does the conjugacy of the $\X$-submaximal subgroups imply the existence of a nilpotent $\pi$-Hall subgroup  for $\pi=\pi(\X)$? A posteriori, with the results of this paper in mind, it is easy to give a counterexample: for $\X=\mathfrak{G}_{\{3,7\}}$, the $\X$-maximal (and $\X$-submaximal) subgroups of $\PSL_2(7)$ are conjugate, while  $\PSL_2(7)$ has no nilpotent $\{3,7\}$-Hall subgroups.  However, in the natural way,  Wielandt's question \cite[offene
Frage, p.~37~(643)]{Wie4} gives rise to the following

\begin{Probl}\label{WieProb2}
{\rm
In what finite simple groups are
all $\X$-(sub)\-ma\-xi\-mal subgroups conjugate?}
\end{Probl}

A complete solution to Problem~\ref{WieProb2} turns out to be a by-product of our approach to Problem~\ref{WieProb1} outlined in~\ref{sub3} and the previous results on $\D_\pi$-groups.

\subsection{ $\D_\X$-groups and $\D_\pi$-groups}\label{sub3}

 We need  the concept of a $\D_\X$-group
introduced in  \cite{GR2} which is natural in the context of our article.

\begin{Def}\label{DX def} A finite group $G$ is called a {\it $\D_\X$-group} (we also say
that $G$  {\it belongs to $\D_\X$} or $G$ {\it satisfies  $\D_\X$}  and write ${G\in\D_\X}$), if ${\k_\X(G)=1}$, i.\,e. the
$\X$-ma\-xi\-mal subgroups of $G$ are conjugate.
\end{Def}

Obviously, Problem~\ref{WieProb1} can be rewritten in the form: {Is it true that an $\X$-submaximal subgroup of a $\D_\X$-group is always $\X$-maximal?} Also, it can be equivalently reformulated as
follows, see   \cite[Theorem~2]{GR2}:

\begin{Probl}\label{WieProb1'} {\rm Is an extension of a $\D_\X$-group by a
$\D_\X$-group always a $\D_\X$-group? In other words, if $N\trianglelefteqslant G$ and $N, G/N\in \D_\X$, is it true that $G\in\D_\X$?
}
\end{Probl}

If ${\X=\mathfrak{G}_\pi}$ is
the class of $\pi$-groups, we write $\D_\pi$ instead of~$\D_\X$. The statement ${G\in \D_\pi}$ is equivalent to the complete analog of the Sylow theorem for $\pi$-subgroup of~$G$ which means the existence and conjugacy of $\pi$-Hall subgroups and the fact that every $\pi$-subgroup is contained in a $\pi$-Hall subgroup. The notation $\D_\pi$ was introduced by P.\,Hall~\cite{Hall}. The particular case of Problem~\ref{WieProb1'} for  $\X=\mathfrak{G}_\pi$ was
first formulated in the one-hour talk by Wielandt at the 13th International Congress of
Mathematicians in Edinburgh in~1958~\cite{Wie1}.\footnote{Moreover, this case of Problem~\ref{WieProb1'} is mentioned in the
surveys \cite{ChuShem,Shem7,Wie5} and in the textbooks
\cite{Suz,ShemBook,GuoBook,GuoBook1}, and was also included by L.\,Shemetkov into the
``Kourovka Notebook'' \cite[Problem 3.62]{Kour}.}
Problem~\ref{WieProb1'} in this case
for  ${\X=\mathfrak{G}_\pi}$
is solved in the affirmative using the classification of finite simple groups (see
\cite[Theorem~6.6]{Surv}). We are going to use this result to solve Problems~\ref{WieProb1} and~\ref{WieProb1'} in general situation. Moreover, for every $\pi$, the simple $\D_\pi$-groups are classified \cite[Theorem~6.9]{Surv}. This classification allows us to solve Problem~\ref{WieProb2}.

 The possibility to reduce the case of an arbitrary complete class $\X$ to the well-studied case $\X=\mathfrak{G}_\pi$ and thereby to solve Problems~\ref{WieProb1}, \ref{WieProb2}, and~\ref{WieProb1'} and to obtain Theorem~\ref{Reduktionssatz_main} is ensured by the following theorem   which is the main result of the paper.

\begin{Theo}\label{DX_means_Dpi} {\rm (mod CFSG)} Let $\X$ be a complete class of finite groups, let $\pi=\pi(\X)$, and let $G$ be a finite simple group. Then $$G\in\D_\X \text{ if and only if either } G\in\X \text{ or } \pi(G)\nsubseteq \pi \text{ and } G\in\D_\pi.$$ In particular, if $G\in\D_\X$ then $G\in\D_\pi$.
\end{Theo}

Since $G\in\D_\pi$ if $G$ has a nilpotent $\pi$-Hall subgroup \cite{Wie1},  Theorem~\ref{DX_means_Dpi} can be considered as a weakening of the original Wielandt's conjecture that the conjugacy of the $\X$-submaximal subgroups of a simple group $G\notin\X$ implies the existence of a nilpotent $\pi$-Hall subgroup \cite[offene
Frage, p.~37~(643)]{Wie4}. Moreover, it follows from \cite[Lemma~4]{GRV} that if $G$ is a simple $\D_\pi$-group and $\pi(G)\nsubseteq \pi$ then either ${2\notin\pi\cap\pi(G)}$ or ${3\notin\pi\cap\pi(G)}$ and the description of $\pi$-Hall subgroups of finite simple groups \cite[Appendix~1]{Surv} implies that $G$ possesses a solvable $\pi$-Hall subgroup $H$ of Fitting height at most~2 (i.\,e. $H$ is an extension of a nilpotent group by nilpotent one).

Briefly, the proof of Theorem~\ref{DX_means_Dpi} is based on the following simple consequence of the Sylow theorem: {\it in a $\D_\X$-group, every $\X$-maximal subgroup is a $\pi$-Hall subgroup for ${\pi=\pi(\X)}$}. We prove the main part `only if' of the theorem considering case by case the simple groups $G$ with a Hall subgroup $H$ (the Hall subgroups of simple groups are known, see \cite[Appendix~1]{Surv}). We show that if $G\notin \D_\pi$ where $\pi=\pi(H)$ is the set of prime divisors of $|H|$, then $G$ contains a subgroup $U$ such that $U$ is not conjugate to a subgroup of $H$, while every composition factor of $U$ is isomorphic to a section (i.\,e. a quotient of a subgroup) of~$H$. It means that if $H\in\m_{\X}(G)$ then $U\in \X$, while there are no $\X$-maximal subgroups of $G$ containing $U$ and conjugate to $H$. In particular, $G\notin\D_\X$. This approach explains the rather large size of this article.

 In~\cite{GR2}, Problems~\ref{WieProb1} and~\ref{WieProb1'} are reduced to simple groups. Using \cite[Theorems~1 and~2]{GR2} and the mentioned above results on $\D_\pi$-groups, we obtain from Theorem~\ref{DX_means_Dpi} the following statements which solve, respectively, the equivalent Problems~\ref{WieProb1'} and~\ref{WieProb1}.

\begin{Cor}\label{DX_Ext}  {\rm (mod CFSG)}  Let $\X$ be a complete class. Assume that $N$ is a normal subgroup of a finite group $G$. Then $G\in\D_\X$ if and only if $N\in\D_\X$ and $G/N\in\D_\X$.
\end{Cor}

\begin{Cor}\label{ConjMax=ConjSubmax} {\rm (mod CFSG)}  Let $\X$ be a
complete class. Then the conjugacy of the
$\X$-ma\-xi\-mal subgroups of a finite group is equivalent to the
conjugacy of the $\X$-sub\-ma\-xi\-mal subgroups.
\end{Cor}

Corollary~\ref{DX_Ext} means that a finite group is a $\D_\X$-group if and only if its every  composition factor is a $\D_\X$-group. In order to describe all $\D_\X$-groups for given $\X$ it is sufficient to list the simple $\D_\X$-groups. This is done in the following statement which also solves Problem~\ref{WieProb2}.

\begin{Cor}\label{ClassificationDX} {\rm (mod CFSG)}  Let $\X$ be a complete class of finite groups. Then, for a finite simple group~$S$, the following statements are equivalent:

 \begin{itemize}
   \item[$(i)$] All $\X$-sub\-ma\-xi\-mal subgroups are conjugate;
   \item[$(ii)$] All $\X$-ma\-xi\-mal subgroups are conjugate $($i.e. $S\in\D_\X)$;
   \item[$(iii)$] Either $S\in \X$, or the pair $(S,\X)$ satisfies one of Conditions {\rm I--VII} in Appendix.
 \end{itemize}
 \end{Cor}

 Conditions I--VII are arithmetic and given in terms of natural parameters of simple groups. These conditions appear in \cite{DpiCl,R2,R4,R5} as necessary and sufficient for a simple group~$S$ to satisfy~$\D_\pi$, see also \cite[Theorem~6.9
and Appendix~2]{Surv}. Note that Condition~I here differs from Condition~I
in~\cite{R4,Surv}. In these articles, Condition~I includes the case
$\pi(S)\subseteq\pi$.

Theorem~\ref{Reduktionssatz_main} immediately follows from Corollaries~\ref{DX_Ext}--\ref{ClassificationDX} and Wielandt's theorem \cite[15.4]{Wie4} which states that the conjugacy of $\X$-submaximal subgroups implies the $\X$-Reduktionssatz.

For ${\pi=\pi(\X)}$, the classes $\D_\X$ and $\D_\pi$ are closely related. Although it is easy to understand that they are different in general (the groups in $\mathfrak{G}_\pi\setminus\X$ are $\D_\pi$-groups but are not  $\D_\X$-groups), it follows by induction from Theorem~\ref{DX_means_Dpi} and Corollary~\ref{DX_Ext} that $$\D_\X\subseteq\D_\pi.$$  This observation allows us to establish rather unobvious properties of $\X$-subgroups in the groups with conjugate $\X$-maximal subgroups. So, in these groups, every $\pi$-subgroup belongs to~$\X$. Applying the results of \cite{Manz,VMR,MRV} we conclude that every subgroup $K$ of  $G\in\D_\X$ which contains an $\X$-maximal subgroup of $G$ is itself a $\D_\X$-group. In particular, an $\X$-maximal subgroup in $K$ is always $\X$-maximal in~$G$.

\subsection{Remarks and open problems}\label{sub4}

 \noindent $1^\circ$ It follows from Corollary~\ref{DX_Ext} that the $\D_\X$-groups, like solvable or nilpotent groups, form a Fitting class. Recall that a {\it Fitting class} is a class of finite groups $\mathfrak{F}$ such that\begin{itemize}
                           \item $N\in \mathfrak{F}$ if $N\trianglelefteqslant G$ and $G\in \mathfrak{F}$ and
                           \item $G\in \mathfrak{F}$ if $G=MN$ where $M,N\trianglelefteqslant G$ and $M,N\in \mathfrak{F}$.
                         \end{itemize}
 For such
 $\mathfrak{F}$, every group $G$ has the {\it $\mathfrak{F}$-radical} $G_{\mathfrak{F}}$, i.\,e. the largest normal $\mathfrak{F}$-subgroup.
 \begin{itemize}
 \item Let $D$ be the $\D_\X$-radical of a group~$G$. Then $D$ coincides with the subgroup generated by all (sub)normal subgroups $U$ of $G$ such that every composition factor of $U$ either is an $\X$-group or satisfies one of Conditions I--VII.
         \end{itemize}
         By Theorem~\ref{Reduktionssatz_main},
 there is a bijection between the conjugacy classes of $\X$-maximal subgroups in $G$ and $G/D$,
   and we can study the members of $\m_\X(G/D)$ instead of $\m_\X(G)$.

  Considering the $\D_\X$-radical $D$ specifically from the point of view of reduction theorems, we can say that $D$ is an absolute radical of $G$ in the sense that $D$ is defined by its inner structure and does not depend on the way of embedding in~$G$.  {\it Is there a subgroup of a group $G$ which is the largest among the normal subgroups $N$ such that ${\k_\X(G)=\k_\X(G/N)}$}? Equivalently, {\it is it true that if $M$ and $N$ are normal subgroups of a finite group $G$ then ${\k_\X(G)=\k_\X(G/M)=\k_\X(G/N)}$ always implies  ${\k_\X(G)=\k_\X(G/MN)}$}? Obviously, an affirmative answer to this question would follow from the validity of the following conjecture: {\it if $N$ is a normal subgroup of a group $G$ then $\k_\X(G)=\k_\X(G/N)$ always implies  ${\k_\X(N)=1}$, i.\,e. ${N\in \D_\X}$}. The authors do not know any counterexamples to this conjecture.

  \medskip

 \noindent $2^\circ$   The concept of an $\X$-submaximal subgroup is useful not only in connection with reduction theorems. In his plenary talk~\cite{Wie3} at the conference on finite groups in Santa-Cruz in 1979, Wielandt suggested a general program on studying $\X$-maximal subgroups of finite groups, and $\X$-submaximal subgroups are  a cornerstone of this program.  We have mentioned that the image under an epimorphism of an $\X$-maximal subgroup is not an $\X$-ma\-xi\-mal subgroup of the image. Also,
 \begin{itemize}
\item the intersection of an  $\X$-ma\-xi\-mal subgroup $H$ with a normal subgroup
$N$ of  $G$ is not an $\X$-ma\-xi\-mal in~$N$, in general.  For example,  a Sylow
2-subgroup $H$ of $G=\PGL_2(7)$ is $\X$-maximal in $G$ for  $\X=\mathfrak{G}_{\{2,3\}}$ but $H\cap
N\notin\m_{\X}(N)$ for $N=\PSL_2(7)$, because $H\cap N$ is a Sylow 2-subgroup of $N$ and is contained in a subgroup of order~$24$
\cite[Example~2, p.~170]{Suz}.
\end{itemize}
In contrast to $\X$-maximal subgroups, $\X$-submaximal ones have the following evident inductive property:
 \begin{itemize}
\item {\it If $H\in \sm_\X(G)$ and $N$ is a (sub)normal subgroup of $G$ then $H\cap N \in \sm_\X(N)$.}
\end{itemize}
 Notice that not every $\X$-subgroup of a group is $\X$-submaximal. An obstruction here is the Wielandt--Hartley theorem. The strong version of this theorem announced in \cite[5.4(a)]{Wie3} and proven in \cite[Theorem~2]{RSV} states:\begin{itemize}
\item {\it If $H\in \sm_\X(G)$ then $N_G(H)/H$ contains  no non-trivial $\X$-subgroups.}
\end{itemize}
For the case where $H$ is strongly $\X$-submaximal, this theorem was proven in \cite[13.2]{Wie4} and \cite[Lemmas 2 and~3]{Hart}.
It is due to the Wielandt--Hartley theorem and the above-mentioned inductive property that the concept of an $\mathfrak{X}$-submaximal subgroup becomes useful and efficient. For example, it helps to easily see that the $\mathfrak{X}$-maximal subgroups are determined uniquely up to conjugacy by their projections on the sections of a subnormal series \cite[5.4(c)]{Wie3}, \cite[Corollary~1]{RSV}. It was shown in \cite{GR1,GR_surv} that the knowledge of $\mathfrak{X}$-submaximal subgroups in simple groups for a given class $\mathfrak{X}$ would make it possible to inductively construct the $\mathfrak{X}$-maximal subgroups in an arbitrary finite group and, consequently, would make great progress in solving the general problem of determining the $\X$-maximal subgroups. Thus, the central problem in Wielandt's program and in the topic related to the search of $\mathfrak{X}$-maximal subgroups for a given complete class~$\mathfrak{X}$ is the description of $\mathfrak{X}$-submaximal subgroups in simple groups which is equivalent to the description of  $\mathfrak{X}$-maximal subgroups in the automorphism groups of simple groups \cite[5.3]{Wie3}.
As a part of this problem pointed out by Wielandt himself \cite[Frage~g]{Wie3}, the classification of $\X$-sub\-ma\-xi\-mal subgroups
in minimal non-solvable groups is completed in \cite{GR_SubmaxMinNonSolv}.

  \medskip

 \noindent $3^\circ$
 In light of Theorem~\ref{Reduktionssatz_main} and the importance of the $\X$-submaximal subgroups, it is natural to ask about the validity of reduction theorems for these subgroups themselves. Wielandt thought that such theorems hold, which can be seen from his latest diary entries. For example, \cite[p.~57, Satz~1]{Wie_diary} states for $\X=\mathfrak{G}_\pi$: if a normal subgroup $N$ of a group~$G$ is such that the $\X$-submaximal subgroups are conjugate in every composition factor of~$N$, then
there is a natural one-to-one correspondence between the $G$-  and $\overline{G}$-conjugacy classes of  $\m_\X(G)$ and  $\m_\X(\overline{G})$, where $\overline{G}=G/N$, as well as between the classes of  $\sm_\X(G)$ and  $\sm_\X(\overline{G})$.
However,  the
`as well as' part is incorrect even in the case when $N$ is a $\pi$-se\-par\-ab\-le group. Counterexamples are given in \cite{RZ1,RZ2}. The number of conjugacy classes of $\X$-submaximal subgroups in $G$ can be both strictly greater and strictly smaller than in~$G/N$.
The following statement is a consequence of our results. It can be considered as a weak analog of Theorem~\ref{Reduktionssatz_main} for $\X$-submaximal subgroups.

\begin{Cor}\label{moduloRadicals} {\rm (mod CFSG)}
Let $\X$ be a complete class of finite groups, let $\mathfrak{F}$ be a Fitting class such that $\mathfrak{F}\subseteq\D_\X$, and let $G$ be a finite group. Denote by $\overline{\phantom{x}}:G\rightarrow G/G_{\mathfrak{F}}$ the canonical epimorphism. If $H\in\sm_\X(G)$ then $\overline{H}\in \sm_\X(\overline{G})$. In particular, if $k$ and $\overline{k}$ are the numbers of conjugacy classes of $\X$-submaximal subgroups in $G$ and $\overline{G}$, respectively, then
$\overline{k}\ge k.$
\end{Cor}

The classes of nilpotent groups, solvable groups, $\X$-groups, groups without non-trivial $\X$-subgroups, $\X$-separable groups, $\D_\X$-groups are examples of Fitting classes  $\mathfrak{F}$ satisfying the assumption of Corollary~\ref{moduloRadicals}. In some sense, $\X$-submaximal subgroups can be studied modulo  the $\mathfrak{F}$-radical for such  $\mathfrak{F}$. Of course, Corollary~\ref{moduloRadicals} does not look as impressive as Theorem~\ref{Reduktionssatz_main}. Examples  from \cite{RZ1,RZ2} show that we can not replace the $\mathfrak{F}$-radical $G_{\mathfrak{F}}$ in the assumption of  Corollary~\ref{moduloRadicals} with an arbitrary normal $\D_\X$-subgroup~$N$, because the image in $G/N$ of an $\X$-submaximal subgroup of $G$ is not $\X$-submaximal in $G/N$ in general.  However, we do not know of any such examples where~$N$ is characteristic. Notice that even under the assumption of Corollary~\ref{moduloRadicals}, $G/G_{\mathfrak{F}}$ may possesses an $\X$-submaximal subgroup which is not the image of  any $\X$-submaximal subgroup of~$G$ and the inequality in Corollary~\ref{moduloRadicals} may be strict.

\section{Preliminaries}

\subsection{ Notation}

\noindent According to \cite{Bray,Atlas, KL}, we use the following notation.

\begin{itemize}
\item[] $\varepsilon$ and $\eta$ always  denote either $+1$ or $-1$ and the
    sign of this number. Sometimes (in the notation of orthogonal groups of
    odd dimension) $\eta$ can be used as an empty symbol.

  \item[] $n$ denotes the cyclic group of order $n$, where $n$ is a
      positive integer.

   \item[] $A^n$  denotes the direct product of $n$ copies of $A$. In
       particular,

   \item[] $p^n$  denotes the elementary abelian group of order $p^n$,
       where $p$ is a prime.

%
\item[] $\Sym(\Omega)$ denotes the symmetric group on $\Omega$.

\item[] $\Sym_n$ is the symmetric group of degree $n$, i.e.
    $\Sym_n=\Sym(\Omega)$, where $\Omega=\{1,2,\dots,n\}$.

\item[] $\Alt_n$  denotes the alternating group of degree $n$.

\item[] $\GL_n(q)$ or $\GL^+_n(q)$ denotes the general linear group of
    degree $n$ over a field of order~$q$.

\item[] $\SL_n(q)$   or $\SL^+_n(q)$ denotes the special linear group of
    degree $n$ over a field of order~$q$.

\item[] $\PSL_n(q)$  or $\PSL^+_n(q)$ denotes the projective special linear
    group of degree $n$ over a field of order~$q$.

\item[] $\PGL_n(q)$  or $\PGL^+_n(q)$ denotes the projective general linear
    group of degree $n$ over a field of order~$q$.

\item[] $\GU_n(q)$  or $\GL^-_n(q)$ denotes the general unitary group of
    degree $n$ over a field of order~$q$.

\item[] $\SU_n(q)$  or $\SL^-_n(q)$ denotes the special unitary group of
    degree $n$ over a field of order~$q$.

\item[] $\PSU_n(q)$  or $\PSL^-_n(q)$ denotes the projective special
    unitary group of degree $n$ over a field of order~$q$.

\item[] $\PGU_n(q)$  or $\PGL^-_n(q)$ denotes the projective general
    unitary group of degree $n$ over a field of order~$q$.

\item[] $\mathrm{SO}_n^\varepsilon(q)$ is the orthogonal group of degree
    $n$ over a field of order $q$, where $\varepsilon\in\{+1,-1\}$ for $n$
    even and $\varepsilon$ is an empty symbol for $n$ odd.

\item[] $\Omega_n^\varepsilon(q)$ is the derived subgroup of
    $\mathrm{SO}_n^\varepsilon(q)$.

\item[] $\mathrm{P}\Omega_n^\varepsilon(q)$ is the reduction of
    $\Omega_n^\varepsilon(q)$  modulo scalars.

\item[] $\Sp_n(q)$  denotes the symplectic group of even degree $n$  over a
    field of order~$q$.

\item[] $\PSp_n(q)$  denotes the projective symplectic group of even degree
    $n$  over a field of order~$q$.

\item[] $r^{1+2n}$  denotes an extra special group of  order $r^{1+2n}$,
    where $r$ is a prime.

\item[]$A:B$ means a split extension of a group $A$ by a group $B$ ($A$ is
    normal).

\item[]$A\nonsplitext  B$ means a non-split extension of a group $A$ by a group $B$ ($A$ is normal).

\item[]$A\arbitraryext  B$ means an arbitrary (split or non-split extension) of a group $A$ by a group $B$ ($A$ is normal).

  \item[]$A^{m+n}$ means $A^m:A^{n}$.

  \item[] ${\rm P}G$, for a linear group $G$, means the reduction of $G$ modulo scalars.

\item[] $\X$ is a complete class of groups.

\item[] $\mathfrak{S} $ is a class of all solvable groups.

\item[] $\D_\X$ is a class of groups with all maximal $\X$-sub\-gro\-ups
    conjugate.

\item[] $\pi$ is a set of primes.

\item[] $\mathfrak{S}_\pi$ is a class of all solvable $\pi$-gro\-ups

\item[] $\G_\pi$ is a class of all $\pi$-gro\-ups.

\item[] $\D_\pi$ is a class of groups with all maximal $\pi$-sub\-gro\-ups
    conjugate, i.e. $\D_\pi=\D_{\G_\pi}$.

\item[] $\E_\pi$ is a class of groups possessing $\pi$-Hall subgroups, i.e. $G\in \E_\pi$, if $\Hall_\pi(G)$ is nonempty.

 \item[] $G_\X$ means the $\X$-radical of $G$, i.e. the subgroup generated
     by all normal $\X$-sub\-gro\-ups of $G$. In particular,

  \item[] $G_{\mathfrak{S}}$ means the solvable radical of $G$.

  \item[] $O_\pi(G)$ means the $\pi$-radical of $G$, the subgroup
      generated by all normal $\pi$-subgroups of $G$, i.e.
      $O_\pi(G)=G_{\G_\pi}$.

\item[] $\mu(G)$ denotes the degree of the minimal faithful permutation
    representation of a finite group~$G$, i.e. the smallest $n$ such that
    $G$ is isomorphic to a subgroup of~$\Sym_n$.

\item[] $\X$-Hall subgroup, is a subgroup $H$ of $G$ such that $H$ is an
    $\X$-subgroup and a $\pi(\X)$-Hall subgroup.

\item[] $\Hall_\X(G)$ is the set of all $\X$-Hall subgroups of $G$, i.e.
    $\Hall_\X(G)=\X\cap \Hall_{\pi(\X)}(G)$.

\end{itemize}

\subsection{ Known properties of $\pi$-Hall subgroups, $\D_\pi$- and $\D_\X$-groups}

\begin{Lemma}\label{HallSubgroup} {\em \cite[Lemma~1]{Hall}}
Let $N$ be a normal subgroup and let $H$ be a $\pi$-Hall subgroup of a group $G$. Then
${H\cap N\in \Hall_\pi(N)}$ and  $HN/N\in \Hall_\pi(G/N).$
\end{Lemma}

\begin{Lemma}\label{Gross}
{\em \cite[Theorem~A]{Gross3}} If $2\notin\pi$ and $G\in\E_\pi$, then every two $\pi$-Hall subgroups of $G$ are conjugate.
\end{Lemma}

\begin{Lemma}\label{DpiExt}
{\em \cite[Theorem~7.7]{Hall3'}, \cite[Theorem~6.6]{Surv}} Let $N$  be a
normal subgroup of  $G$. Then  $G\in\D_\pi$ if and only if  $N\in \D_\pi$ and
$G/N\in \D_\pi$.
\end{Lemma}

\begin{Lemma}\label{DpiCrit} {\em \cite[Theorem 3]{R4}} Let $\pi$ be a set of primes and $S$ a simple group. Then $S\in \D_\pi$ if and only  if either $S$ is a $\pi$-group or  $(S,\X)$ for $\X=\mathfrak{G}_\pi$ satisfies one of Conditions~{\rm I--VII}.
\end{Lemma}

\begin{Lemma}\label{HallXSubgroups} {\em \cite[Proposition~1]{GR2}} Let $\X$ be a complete class, $\pi=\pi(\X)$, and $G\in\D_\X$. Then
$
\m_\X(G)=\Hall_\X(G)\subseteq\Hall_\pi(G).
$
In particular, $G\in\E_\pi$.
\end{Lemma}

\begin{Lemma}\label{DX_Normal_and_Quot} {\rm \cite[Theorem~1]{GR2}} If $G\in\D_\X$ and $N\trianglelefteqslant G$ then $N\in\D_\X$ and $G/N\in\D_\X$.
 \end{Lemma}

 \begin{Lemma}\label{DXExtEquiv} {\rm \cite[Theorem~2]{GR2}} For  complete  $\X$, the following statements are equivalent.
\begin{itemize}
 \item[$(1)$]  The elements of ${\rm sm}_\X(G)$ are conjugate in any  $G\in\D_\X$.
  \item[$(2)$]  $\sm_\X(G)=\m_\X(G)$  for any  $G\in\D_\X$.
  \item[$(3)$]   $\D_\X$ is closed under taking extensions.
  \item[$(4)$]  $\Aut S\in\D_\X$ for every simple  $S\in\D_\X$.
  \item[$(5)$]  The elements of ${\rm sm}_\X(S)$ are conjugate in any simple  $S\in\D_\X$.
\end{itemize}
 \end{Lemma}


\begin{Lemma}\label{Bijection}
{\rm \cite[12.9]{Wie4}} {\it Let $N$ be a normal $\X$-separable subgroup of
$G$. Then the map given by the rule $M\mapsto MN/N$ is a surjection between
sets $\m_\X(G)$ and $\m_\X(G/N)$. Moreover, this map induces a bijection
between the sets of conjugacy classes of $\X$-ma\-xi\-mal subgroups of $G$
and $G/N$. In particular, $G\in \D_\X$ if and only if $G/N\in \D_\X$.}
\end{Lemma}


\begin{Lemma}\label{Ext} {\rm \cite[Lemma~2.1(e)]{NumbCl}} {\it
 Let $N$ be a normal subgroup of $G$ and $\pi(G/N)\subseteq \pi$. Assume $N$ contains a $\pi$-Hall subgroup $H_0$. Then the following statements are equivalent.
 \begin{itemize}
  \item[$(1)$] There is $H\in\Hall_\pi(G)$ such that $H_0=H\cap N$.
  \item[$(2)$] For every $g\in G$ there exists $x\in N$ such that $H_0^g=H_0^x$.
 \end{itemize}}
\end{Lemma}

\begin{Lemma}{\rm (The Wielandt--Hartley theorem, strong form) \cite[5.4(a)]{Wie3}, \cite[Theorem 2]{RSV}}\label{WH-strong}
  Let  $K$ be an $\X$-submaximal subgroup of a finite group $G$. Then $N_G(K)/K$
  is a $\pi(\X)'$-group.
\end{Lemma}


\subsection{Arithmetic Lemmas}

For an odd integer $q$, denote by $\varepsilon(q)$ the number
$\varepsilon=\pm1$ such that $q\equiv \varepsilon\pmod 4$.

If $r$~is an odd prime and  $k$~is an integer not divisible by $r$, then
$e(k,r)$ is the smallest
 positive integer  $e$ with $k^e \equiv 1 \pmod r$.
 So, $e(k,r)$ is the multiplicative order of $k$ modulo~$r$.

For a natural number $e$ set
 $$e^*=\left\{
\begin{array}{ll}
2e & \text{ if }  e\equiv 1\pmod 2,\\
e & \mbox{ if   }  e\equiv 0\pmod 4,\\
e/2 & \mbox{ if  }  e\equiv 2\pmod 4.
\end{array}
\right.
$$
It follows from  the definition that if $e$ divides an even number $n$ then
$e^*$ divides $n$ again. Moreover,  $e^{**}=e$
for every~$e$.

For a real $x$, the integer part $[x]$ of $x$ is defined as a
unique integer $k$ such that $$k\le x< k+1.$$

The following lemma is evident.

\begin{Lemma}\label{intpart}
If $m$ is a positive integer and $x$ is a real then
 \begin{equation*}
   [{[x]}/{m}]=[{x}/{m}].
 \end{equation*}
\end{Lemma}

The next result may be found in \cite{Weir1955}.

\begin{Lemma}\label{r4ast'}{\em (\cite{Weir1955}, \cite[Lemmas~2.4 and
2.5]{Gross2})}
 Let $r$ be an odd prime, $k$~an integer not divisible by $r$,  and  $m$~a
positive integer. Denote $e(k,r)$ by $e$.

 Then the following identities hold.

$$(k^m-1)_r=\left\{
\begin{array}{ll}
(k^e-1)_r (m/e)_{r} & \text{ if } e \text{ divides } m,\\
1 & \text{ if } e \text{ does not divide } m;\end{array}\right.$$

$$(k^m-(-1)^m)_r=\left\{
\begin{array}{ll}
(k^{e^*}-(-1)^{e^*})_r (m/{e^*})_{r} & \text{ if } e^* \text{ divides } m,\\
1 & \text{ if }  e^* \text{ does not divide } m.\end{array}\right.$$

%

\end{Lemma}

 \begin{Lemma}\label{arithm} Let $q>1$ and $n$~be positive integers,  let $r$~be an odd prime such that $(q,r)=1$, and let $e=e(r,q)$. Then the following statements hold{\rm:}

\begin{itemize}
  \item[{\rm(1)}] $(n!)_r=r^\alpha$, where $\alpha =\displaystyle\sum\limits_{i=1}^\infty [n/r^i]${\rm;}

\item[{\rm(2)}]  $\displaystyle\prod\limits_{i=1}^n (q^i-1)_r=(q^e-1)_r^{[n/e]}([n/e]!)_r${\rm;}

\item[{\rm (3)}] $\displaystyle\prod\limits_{i=1}^{m}(k^i-(-1)^i)_r=(k^{e^*}-(-1)^{e^*})_r^{[m/e^*]}([m/e^*]!)_r${\rm;}

\item[{\rm(4)}]  $\displaystyle\prod\limits_{i=1}^n (q^i-1)_r=(n!)_r$ if and only if $e=r-1$, $(q^{r-1}-1)_r=r$ and $[n/r]=[n/(r-1)]${\rm.}

\item[{\rm(5)}]  $\displaystyle\prod\limits_{i=1}^{m}(q^i-(-1)^i)_r=(n!)_r$
    if and only if $e^*=r-1$, $(q^{(r-1)^*}-(-1)^{(r-1)^*})_r=r$ and
    $[n/r]=[n/(r-1)]${\rm.}
\end{itemize}
\end{Lemma}

 \noindent{\sc Proof}.
The statement (1) is well-known (see, for example, \cite[Lemma~2]{Thomp}). Statements (2) and (3) follow from Lemma~\ref{r4ast'}.

Now we prove (4). Let $A=\displaystyle\prod\limits_{i=1}^n (q^i-1)_r$. Then
by (2) and in view of the Little Fermat  Theorem,
 \begin{multline}\label{ne}
   A=(q^e-1)_r^{[n/e]}([n/e]!)_r\geq r^{[n/e]}([n/e]!)_r\geq\\ r^{[n/(r-1)]}([n/(r-1)]!)_r\geq r^{[n/r]}([n/r]!)_r=  r^\beta,
 \end{multline}
where by (1) and, in view of Lemma~\ref{intpart} for $x=n/r$ and $m=r^i$, we have
$$
\beta=[n/r]+\sum\limits_{i=1}^\infty \left[{[n/r]}/{r^i}\right]=[n/r]+\sum\limits_{i=1}^\infty \left[{n}/{r^{i+1}}\right]=\sum\limits_{i=1}^\infty \left[{n}/{r^{i}}\right]=\log_r(n!)_r.
$$
Therefore, $ A\geq (n!)_r. $ Moreover, this inequality becomes an equality if
and only if all inequalities in~(\ref{ne}) are equalities, i.e. if and only
if $$ r-1= e,\,\,\, (q^{e}-1)_r= r,\,\,\,\text{ and }\,\,\, [n/(r-1)]= [n/r].
$$
This implies (4).

Now we prove (5). Let $A'=\displaystyle\prod\limits_{i=1}^n (q^i-(-1)^i)_r$. Since $r$ is odd in view of the Little Fermat Theorem  $e$ divides the even number $r-1$. Consequently, $e^*$ also divides $r-1$   and,  by (3),
\begin{multline}\label{ne'}
   A'=(q^{e^*}-(-1)^i)_r^{[n/e^*]}([n/e^*]!)_r\geq r^{[n/e^*]}([n/e^*]!)_r\geq\\
    r^{[n/(r-1)]}([n/(r-1)]!)_r\geq r^{[n/r]}([n/r]!)_r=  r^\beta,
 \end{multline}
where $\beta$ is as above.

Therefore, $ A'\geq (n!)_r. $ Again this inequality becomes an equality if
and only if all the inequalities in~(\ref{ne'}) are equalities, i.e. if and
only if one of the equalities $$r-1= e^*,\,\,\, (q^{e^*}-(-1)^{e^*})_r=
r,\,\,\,\text{ and }\,\,\, [n/(r-1)]= [n/r].
$$
This implies~(v).
\qed\medskip

\subsection{On Hall subgroups of finite simple groups}

\begin{Lemma}\label{HallSymmetric} {\rm \cite[Theorem~A4 and the notices after it]{Hall}, \cite[Main result]{Thomp}, \cite[Theorem~8.1]{Surv}}
Suppose that $n\geq 5$ and $\pi$ is a set of primes with $|\pi\cap \pi(n!)|>1$ and $\pi(n!)\nsubseteq\pi$. Then
\begin{itemize}
 \item[{\em (1)}] The complete list of possibilities for $\Sym_n$
     containing a $\pi$-Hall subgroup $H$ is given in
     Table~{\em\ref{Symmetric}}.
 \item[{\em (2)}] $M\in\Hall_\pi(\Alt_n)$ if and only if  ${M=H\cap
\Alt_n}$ for some ${H\in\Hall_\pi(\Sym_n)}$.
\end{itemize}
In particular, every proper nonsolvable $\pi$-Hall subgroup of a symmetric group of degree $n$ is isomorphic to a symmetric group of degree $n$ or $n-1$ and has a unique nonabelian composition factor isomorphic to an alternating group of the same degree.
\end{Lemma}

\begin{longtable}{|c|c|c|}\caption{\label{Symmetric}}\\
\hline
$n$&$\pi$&$H\in\Hall_\pi(\Sym_n)$\\ \hline\hline
prime&$\pi((n-1)!)$&$\Sym_{n-1}$\\
$7$&$\{2,3\}$&$\Sym_3\times\Sym_4$\\
$8$&$\{2,3\}$&$\Sym_4\wr\Sym_2$\\ \hline
\end{longtable}


\begin{Lemma}\label{An} {\rm \cite[Proposition~3]{GR2}} {\it Let  $\pi=\pi(\X)$.  Then for $G=\Alt_n$ the following conditions are equivalent.
 \begin{itemize}
   \item[$(1)$] $G\in\D_\X$.
   \item[$(2)$] $G\in\D_\X\cap\D_\pi$.
   \item[$(3)$] either $|\pi\cap\pi(G)|\leq 1$ or $G\in\X$.
   \item[$(4)$] All submaximal $\X$-subgroups are conjugate in~$G$.
 \end{itemize}
 }
\end{Lemma}

\begin{Lemma}\label{HallSpor23} {\em \cite[Theorem 4.1]{DpiCl}, \cite[Theorem~8.2]{Surv}}
Let $G$ be either one of
the $26$ sporadic groups or a Tits group,  $\pi$ be
such that $2\in\pi$, $\pi(G)\nsubseteq\pi$, and $|\pi\cap\pi(G)|>1$, and $H$ be a $\pi$-Hall subgroup of $G$. Then the corresponding intersections $\pi\cap \pi(G)$ and the structure of $H$  are indicated in Table~{\em\ref{tb0}}.
\end{Lemma}

\begin{longtable}{|l|l|r|}\caption{}\label{tb0}\\  \hline
$G$ & $\pi\cap\pi(G)$ &
Structure of $H$ \\
\hline\hline
$M_{11}$  & $\{2,3\}$ & $3^2\splitext Q_8\arbitraryext 2$ \\
          & $\{2,3,5\}$       &  $\Alt_6\arbitraryext 2$\\ \hline
$M_{22}$  & $\{2,3,5\}$ & $2^4\splitext \Alt_6$\\ \hline
$M_{23}$  & $\{2,3\}$         & $2^4\splitext(3\times A_4)\splitext2$\\
& $\{2,3,5\}$ &  $2^4\splitext \Alt_6$\\ & $\{2,3,5\}$       &
$2^4\splitext (3\times \Alt_5)\splitext 2$\\
& $\{2,3,5,7\}$     &  ${\rm           L}_3(4)\splitext 2_2$\\
& $\{2,3,5,7\}$     &  $2^4\splitext \Alt_7$\\
&           $\{2,3,5,7,11\}$  & $M_{22}$\\ \hline
$M_{24}$  & $\{2,3,5\}$    & $2^6 \splitext 3\nonsplitext \Sym_6$\\ \hline
$J_1$     & $\{2,3\}$         & $2\times \Alt_4$\\
 & $\{2,3,5\}$ &           $2\times \Alt_5$\\
          & $\{2,3,7\}$       &  $2^3\splitext 7\splitext 3$\\
          & $\{2,7\}$       &  $2^3\splitext 7$\\ \hline
$J_4$     & $\{2,3,5\}$       &  $2^{11}\splitext (2^6\splitext 3\nonsplitext \Sym_6)$\\
\hline
\end{longtable}

\begin{Lemma}\label{SLU23dim2} {\rm \cite[Lemma~3.1]{NumbCl}, \cite[Lemma~8.10]{Surv}} Let $\pi$ be a set of primes with $2,3\in\pi$.  Assume that
$G\simeq\SL_2(q)\simeq\SL_2^\eta(q)\simeq\Sp_2(q)$,  where
$q$ is a power of an odd prime $p\not\in\pi$, and $\varepsilon=\varepsilon(q)$.  Recall that for a subgroup $A$ of $G$ we denote by
$\P A$ the reduction modulo scalars. Then the following statements hold.

\begin{itemize}
\item[{\em (A)}] If $G\in \E_\pi$ and $H\in \Hall_\pi(G)$, then one of  the following statements holds.
\begin{itemize}
\item[{\em (a)}]  $\pi\cap\pi(G)\subseteq \pi(q-\varepsilon)$, $\P H$ is a
$\pi$-Hall subgroup in the dihedral subgroup $D_{q-\varepsilon}$ of order $q-\varepsilon$ of $\P G$.
All $\pi$-Hall subgroups of this type are conjugate in~$G$.

\item[{\em (b)}] $\pi\cap\pi(G)=\{2,3\}$, $(q^2-1)_{\mbox{}_{\{2,3\}}}=24$, $\P H\simeq \Alt_4$.
All  $\pi$-Hall subgroups of this type are conjugate in~$G$.

\item[{\em (c)}] $\pi\cap\pi(G)=\{2,3\}$, $(q^2-1)_{\mbox{}_{\{2,3\}}}=48$, $\P H\simeq \Sym_4$.
There exist exactly two classes of conjugate subgroups of this type, and
$\PGL_2^\eta(q)$ interchanges  these classes.

\item[{\em (d)}] $\pi\cap\pi(G)=\{2,3,5\}$, $(q^2-1)_{\mbox{}_{\{2,3,5\}}}=120$, $\P H\simeq
\Alt_5$. There exist exactly two classes of conjugate subgroups of this type, and
$\PGL_2^\eta(q)$ interchanges  these classes.
\end{itemize}
\item[{\em (B)}] Conversely, if $\pi$ and $(q^2-1)_\pi$ satisfy one of statements {\em (a)--(d)}, then~${G\in \E_\pi}$.
\item[{\em (C)}] Every $\pi$-Hall subgroup of $\P G$ can be obtained as $\P H$ for some~${H\in \Hall_\pi(G)}$. Conversely, if $\P H\in
\Hall_\pi(\P G)$ and $H$ is a full preimage of $\P H$ in $G$, then~$H\in \Hall_\pi(G)$.
\end{itemize}
\end{Lemma}

\begin{Lemma}\label{GLU23dim2}  {\rm \cite[Lemma~3.2]{NumbCl}, \cite[Corollary~8.11]{Surv}}
Let $G=\GL^\eta_2(q)$, $\P G=G/Z(G)=\PGL_2^\eta(q)$, where $q$ is a power of a
prime $p$, and $\varepsilon=\varepsilon(q)$.  Let $\pi$ be a set of primes such that $2,3\in\pi$ and $p\not\in\pi$. A
subgroup  $H$ of $G$ is a $\pi$-Hall subgroup if and only if $H\cap \SL_2(q)$ is
a $\pi$-Hall subgroup of $\SL_2(q)$, $\vert H:H\cap\SL_2(q)\vert_\pi=(q-\eta)_\pi$, and either statement {\em (a)}, or
statement {\em (b)} of Lemma {\rm
\ref{SLU23dim2}} holds. More precisely, one of the following statements holds.
\begin{itemize}
\item[{\em (a)}]  $\pi\cap\pi(G)\subseteq \pi(q-\varepsilon)$, where
$\varepsilon=\varepsilon(q)$,
$\P H$ is a $\pi$-Hall subgroup in the dihedral group $D_{2(q-\varepsilon)}$ of order
$2(q-\varepsilon)$ of~$\P G$. All $\pi$-Hall subgroups of this type are
conjugate in~$G$.

\item[{\em (b)}] $\pi\cap\pi(G)=\{2,3\}$, $(q^2-1)_{\mbox{}_{\{2,3\}}}=24$, $\P H\simeq \Sym_4$.
All $\pi$-Hall subgroups of this type are conjugate in~$G$.
\end{itemize}
\end{Lemma}

\begin{Lemma}\label{3.3.12} {\rm \cite[Theorem~3.2]{Gross1}, \cite[Theorem~3.1]{Gross4},\cite[Theorem~1.2]{HallChev}, \cite[Theorems~8.3--8.7]{Surv}}
Let $\pi$ be a set of primes and $G$  a group of Lie type over the field $\mathbb{F}_q$ of characteristic~${p\in\pi}$. Assume,  $G\in \E_{\pi}$ and  $H\in\Hall_\pi(G)$. Then one of the following statements holds.

\begin{itemize}
\item[{\rm (1)}] $H=G$.

\item[{\rm (2)}]
$\pi\cap\pi(G)\subseteq\pi(q-1)\cup\{p\}$, $H$ is contained in a Borel subgroup of~$G$ (in particular, $H$ is solvable) and any prime in $\pi\setminus\{p\}$ does not divide the order of the Weyl group of~$G$.

\item[{\rm (3)}] $p=2$, $G=D_l(q)$, the Dynkin diagram of the fundamental root system $\Pi^1$ of $G$ is on Pic.~$1$,  $l$ is a Fermat prime, $(l,q-1)=1$ and $H$ is conjugate to the canonic parabolic maximal subgroup corresponding to the set $\Pi\setminus\{r_1\}$ of fundamental roots.
\begin{figure}[b]
\begin{center}
\unitlength=8mm \begin{picture}(9,4)
   \put(0,2){\circle* {0.2}}   \put(0,2.5){\makebox(0,0){$r_1$}}
   \put(2,2){\circle* {0.2}}   \put(2,2.5){\makebox(0,0){$r_2$}}
   \put(4,2){\circle* {0.2}}   \put(4,2.5){\makebox(0,0){$r_{l-3}$}}
   \put(6,2){\circle* {0.2}}
 \put(5.8,2.5){\makebox(0,0){$r_{l-2}$}} \put(8,4){\circle*
{0.2}} \put(8.62,4){\makebox(0,0){$r_{l-1}$}} \put(8,0){\circle*
{0.2}} \put(8.5,0){\makebox(0,0){$r_l$}}
 \put(0.1,2){\line(1,0){1.8}}
\put(2.1,2){\line(1,0){0.5}}
\put(3.4,2){\line(1,0){0.5}}
\put(3,2){\makebox(0,0){$\cdots$}}
\put(4.1,2){\line(1,0){1.8}}
\put(6.08,2.08){\line(1,1){1.83}}
\put(6.08,1.92){\line(1,-1){1.83}}
\end{picture}
\centerline{Pic. 1. Dynkin diagram of the root system of  $D_l(q)$.}
\end{center}
\end{figure}

\item[{\rm (4)}] $p=2$, $G={^2D}_l(q)$, the Dynkin diagram of the fundamental root system $\Pi$ of $G$ is on Pic.~$2$,
 $l-1$ is a Mersenne prime, $(l-1,q-1)=1$ and $H$ is conjugate to the canonic parabolic maximal subgroup corresponding to the set $\Pi^1\setminus\{r^1_1\}$ of fundamental roots;
 \begin{figure} \begin{center} \unitlength=9mm
\begin{picture}(9,6.3)
\put(0,4){\circle* {0.2}}
   \put(0,4.5){\makebox(0,0){$r_1$}} \put(2,4){\circle* {0.2}}
   \put(2,4.5){\makebox(0,0){$r_2$}} \put(4,4){\circle* {0.2}}
 \put(4,4.5){\makebox(0,0){$r_{l-3}$}} \put(6,4){\circle* {0.2}}
   \put(6,4.5){\makebox(0,0){$r_{l-2}$}} \put(8,6){\circle*
   {0.2}}   \put(8.5,6){\makebox(0,0){$r_{l-1}$}}
   \put(8,2){\circle* {0.2}}
   \put(8.5,2){\makebox(0,0){$r_l$}}
 \put(0.1,4){\line(1,0){1.8}}
\put(2.1,4){\line(1,0){0.5}}
\put(3.4,4){\line(1,0){0.5}}
\put(3,4){\makebox(0,0){$\cdots$}}
\put(4.1,4){\line(1,0){1.8}}
\put(6.08,4.08){\line(1,1){1.83}}
\put(6.08,3.92){\line(1,-1){1.83}}
\put(8,4){\vector(0,-1){1.9}}
\put(8,4){\vector(0,1){1.9}}
\put(8.5,4){\oval(2.5,5)}
\put(4,3){\vector(0,-1){1.9}}
   \put(0,0){\circle* {0.2}}   \put(0,0.5){\makebox(0,0){$r_1^1$}}
   \put(2,0){\circle* {0.2}}   \put(2,0.5){\makebox(0,0){$r_2^1$}}
   \put(4,0){\circle* {0.2}}
\put(4,0.5){\makebox(0,0){$r_{l-3}^1$}} \put(6,0){\circle* {0.2}}
 \put(5.8,0.5){\makebox(0,0){$r_{l-2}^1$}} \put(8,0){\circle*
{0.2}} \put(8,0.5){\makebox(0,0){$r_{l-1}^1$}}
 \put(0.1,0){\line(1,0){1.8}}
\put(2.1,0){\line(1,0){0.5}}
\put(3.4,0){\line(1,0){0.5}}
\put(3,0){\makebox(0,0){$\cdots$}}
\put(4.1,0){\line(1,0){1.8}}
\put(6,0.1){\line(1,0){2}}
\put(6,-0.1){\line(1,0){2}}
\end{picture}
\centerline{Pic. 2. Dynkin diagram of the root system of
${^2D_l}(q)$.}
\end{center}
\end{figure}

\item[{\rm (5)}]  $G$ is isomorphic to the quotient by the center of ${\rm SL}({ V})$, where ${ V}$~is a vector space of a dimension $n$ over $\mathbb{F}_q$ and $H$ is the image in $G$ under the natural epimorphism of the stabilizer in ${\rm SL}({ V})$ of a series
$$0={ V}_0<{ V}_1<\dots <{ V}_s={ V}$$ such that ${\rm
dim}{ V}_i/{ V}_{i-1}=n_i$, $i=1,2,\dots ,s$, and one of the following conditions holds:
\begin{itemize}
\item[{\rm (a)}]  $n$~is a prime, $(n, q-1)=1$, $s=2$,
$n_1,n_2
\in
\{1,n-1\}$;

\item[{\rm (b)}] $n=4$, $(2\cdot3 ,q-1)=1$, $s=2$,
$n_1=n_2=2$;

\item[{\rm (c)}] $n=5$, $(2\cdot 5 , q-1)=1$, $s=2$,
$n_1,n_2\in\{2,3\}$;

\item[{\rm (d)}] $n=5$, $(2\cdot 3\cdot 5 , q-1)=1$, $s=3$,
$n_1,n_2,n_3\in\{1,2\}$;

\item[{\rm (e)}] $n=7$, $(5\cdot 7 , q-1)=1$, $(3, q+1)=1$, $s=2$,
$n_1,n_2\in\{3,4\}$;

\item[{\rm (f)}] $n=8$, $(2\cdot5\cdot 7 , q-1)=1$, $(3 , q+1)=1$, $s=2$,
$n_1=n_2=4$;

\item[{\rm (g)}] $n=11$, $(2\cdot3\cdot7\cdot 11 , q-1)=1$,
$(5 , q+1)=1$, $s=2$, $n_1,n_2\in\{5,6\}$.
\end{itemize}
\end{itemize}
\end{Lemma}

\begin{Lemma} \label{HallSubgroupsOfLinearAndUnitaryGroups}  {\rm \cite[Lemma~4.3]{NumbCl}, \cite[Theorem~8.12]{Surv}} Assume $G=\SL_n^\eta(q)$ is a special linear or unitary group with the base
field  $\F_q$ of characteristic $p$ and $n\ge2$. Let $\pi$ be a set of primes such that $2,3\in\pi$ and $p\not\in\pi$. Then the following statements
hold.
\begin{itemize}
\item[{\em (A)}] Suppose $G\in \E_\pi$, and $H$ is a
$\pi$-Hall subgroup
of $G$. Then for $G$, $H$, and $\pi$ one of the following statements holds.
\begin{itemize}
\item[{\em (a)}] $n=2$ and one of the statements {\rm (a)--(d)}  of Lemma~{\rm\ref{SLU23dim2}} holds.

\item[{\em (b)}] either $q\equiv \eta \pmod {12}$, or $n=3$ and
 $q\equiv \eta \pmod 4$, $\Sym_n$ satisfies $\E_\pi$, $\pi\cap\pi(G)\subseteq \pi(q-\eta)\cup \pi(n!)$, and if $r\in(\pi\cap\pi(n!))\setminus
\pi(q-\eta)$, then $|G|_r=|\Sym_n|_r$; $H$ is included in  $$M=L\cap G\simeq Z^{n-1}\arbitraryext \Sym_n,$$ where $L=Z\wr
\Sym_n\leq \GL_n^\eta(q)$ and
$Z=\GL_1^\eta(q)$ is a cyclic group of order $q-\eta$. All $\pi$-Hall subgroups of this type are conjugate in~$G$.

\item[{\em (c)}] $n=2m+k$, where $k\in\{0,1\}$, $m\ge 1$,
$q\equiv -\eta \pmod {3}$, $\pi\cap\pi(G)\subseteq \pi(q^2-1)$, the groups $\Sym_m$ and $\GL_2^\eta(q)$ satisfy $\E_\pi$,
and $$M=L\cap
G\simeq(\underbrace{\GL_2^\eta(q)\circ\dots\circ\GL_2^\eta(q)}_{\displaystyle
m \mbox{ \rm times}})\arbitraryext\Sym_m\circ Z,$$ where $L=\GL_2^\eta(q)\wr \Sym_m\times Z\leq \GL_n(q)$ and $Z$ is a cyclic group of order $q-\eta$ for
$k=1$ and $Z$ is trivial for
$k=0$. The subgroup $H$ acting by conjugation on the set of factors in the central product
\begin{equation}\label{prodgl2}
\underbrace{\GL_2^\eta(q)\circ\dots\circ\GL_2^\eta(q)}_{\displaystyle
m \mbox{ \rm times}}
\end{equation} has at most two orbits.
The intersection of $H$ with each factor $\GL_2^\eta(q)$ in \eqref{prodgl2} is a  $\pi$-Hall subgroup in $\GL_2^\eta(q)$. The intersections with
the factors from the same orbit all satisfy the same statement {\em (a)} or {\em (b)} of Lemma~{\rm \ref{GLU23dim2}}. Two $\pi$-Hall subgroups of $M$ are
conjugate in $G$ if
and only if they are conjugate in $M$. Moreover $M$ possesses one, two, or four classes of conjugate $\pi$-Hall subgroups, while all
subgroups $M$
are conjugate in~$G$.

\item[{\em (d)}] $n=4$, $\pi\cap \pi(G)=\{2,3,5\}$,  $q\equiv5\eta\pmod 8$, $(q+\eta)_3=3$,
$(q^2+1)_5=5$, and $H\simeq 4\arbitraryext 2^4\arbitraryext \Alt_6$. In this case, $G$ possesses exactly two classes of conjugate $\pi$-Hall
subgroups of this
type and $\GL_4^\eta(q)$
interchanges these classes.

\item[{\em (e)}] $n=11$, $\pi\cap \pi(G)=\{2,3\}$, $(q^2-1)_{\mbox{}_{\{2,3\}}}=24$, $q\equiv -\eta\pmod3$, $q\equiv\eta\pmod4$, $H$ is
included in a
subgroup $M=L\cap G,$ where $L$ is a subgroup of $G$ of type $\left((\GL_2^\eta(q)\wr \Sym_4)\perp (\GL_1^\eta(q)\wr \Sym_3)\right)$, and
$$H= \left(\left(\left(Z\circ
2\arbitraryext \Sym_4\right)\wr \Sym_4\right)\times
\left(Z\wr \Sym_3\right)\right)\cap G,$$ where $Z$ is a Sylow $2$-subgroup of a cyclic group of order $q-\eta$. All $\pi$-Hall subgroups of
this type are conjugate in~$G$.
\end{itemize}
\item[{\em (B)}] Conversely, if the conditions on $\pi$, $n$, $\eta$, and $q$ in one of statements  {\em (a)--(e)} are satisfied, then
$G\in \E_\pi$.
\end{itemize}
\end{Lemma}

\begin{Lemma}\label{HallSubgroupsOfSymplecticGroups}  {\rm \cite[Lemma~4.4]{NumbCl}, \cite[Theorem~8.13]{Surv}}
Let $G=\Sp_{2n}(q)$ be a symplectic group over a field $\F_q$ of characteristic~$p$. Assume that  $\pi$ is a set of primes
such that  $2,3\in\pi$ and $p\not\in\pi$. Then the following statements hold.
\begin{itemize}
\item[{\em (A)}] Suppose $G\in \E_\pi$ and $H\in \Hall_\pi(G)$. Then both $\Sym_n$
and $\SL_2(q)$ satisfy $\E_\pi$ and
$\pi\cap\pi(G)\subseteq \pi(q^2-1)$. Moreover, $H$ is a $\pi$-Hall subgroup of

\begin{equation}\label{Symplectic}
  M=\Sp_2(q)\wr \Sym_n\simeq
\big(\underbrace{\SL_2(q)\times\dots\times\SL_2(q)}_{\displaystyle n
\mbox{ \rm times}}\big)\splitext \Sym_n\leq G.
\end{equation}
\item[{\em (B)}] Conversely, if both $\Sym_n$ and $\SL_2(q)$ satisfy $\E_\pi$ and $\pi\cap\pi(G)\subseteq
\pi(q^2-1)$, then $M\in \E_\pi$ and every $\pi$-Hall subgroup $H$ of $M$ is a $\pi$-Hall subgroup of $G$.
\item[{\em (C)}] Two $\pi$-Hall subgroups of $M$ are conjugate in $G$ if and only if they are conjugate in~$M$.
\end{itemize}
\end{Lemma}

\begin{Lemma} \label{HallSubgroupsOfOrthogonalGroupsOfEvenDimension}\label{HallSubgroupsOfOrthogonalGroupsOfOddDimension}  {\rm \cite[Lemma~6.7]{NumbCl}, \cite[Theorem~8.14]{Surv}}  Assume that
$G=\Omega_n^\eta(q)$, $\eta\in\{+,-,\circ\}$, $q$ is a power of a prime $p$, $n\ge 7$, $\varepsilon=\varepsilon(q)$. Let $\pi$ be a set of
primes such that  $2,3\in\pi$, $p\not\in\pi$. Then the following statements hold.

\begin{itemize}
\item[{\em(A)}] If  $G$ possesses a $\pi$-Hall subgroup $H$, then
one of the following statements holds.
\begin{itemize}
\item [{\em (a)}] $n=2m+1$, $\pi\cap\pi(G)\subseteq\pi(q-\varepsilon)$, $q\equiv \varepsilon\pmod {12}$, $\Sym_m\in \E_\pi$, and $H$ is a
$\pi$-Hall subgroup in $M=\big(\Oo_2^\varepsilon(q)\wr \Sym_m\times \Oo_1(q)\big)\cap G.$
All $\pi$-Hall subgroup of this type are conjugate.

\item [{\em (b)}] $n=2m$, $\eta=\varepsilon^m$, $\pi\cap\pi(G)\subseteq\pi(q-\varepsilon)$, $q\equiv \varepsilon\pmod {12}$, $\Sym_m\in
\E_\pi$, and $H$ is a $\pi$-Hall
subgroup in $M=\big(\Oo_2^\varepsilon(q)\wr \Sym_m\big)\cap G.$ All $\pi$-Hall subgroup of this type are conjugate.

\item [{\em (c)}] $n=2m$, $\eta=-\varepsilon^m$, $\pi\cap\pi(G)\subseteq\pi(q-\varepsilon)$, $q\equiv \varepsilon\pmod {12}$,
$\Sym_{m-1}\in \E_\pi$, and $H$ is a
$\pi$-Hall subgroup of $M=\big(\Oo_2^\varepsilon(q)\wr \Sym_{m-1}\times \Oo^{-\varepsilon}_2(q)\big)\cap G.$ All
$\pi$-Hall subgroup of this type are conjugate.

\item [{\em (d)}] $n=11$,  $\pi\cap\pi(G)=\{2,3\}$, $q\equiv \varepsilon\pmod {12}$,  $(q^2-1)_\pi=24$, and $H$ is a $\pi$-Hall subgroup
of $M=\big(\Oo_2^\varepsilon(q)\wr \Sym_{4}\times \Oo_1(q)\wr \Sym_3\big)\cap G.$ All
$\pi$-Hall subgroup of this type are conjugate.

\item [{\em (e)}] $n=12$, $\eta=-$, $\pi\cap\pi(G)=\{2,3\}$, $q\equiv \varepsilon\pmod {12}$, $(q^2-1)_\pi=24$, and $H$ is a $\pi$-Hall
subgroup of $M=\big(\Oo_2^\varepsilon(q)\wr \Sym_{4}\times  \Oo_1(q)\wr \Sym_3\times \Oo_1(q)\big)\cap G$.
There exist precisely two classes of conjugate subgroups of this type in $G$,
and the
automorphism of order $2$ induced by the group of similarities of the natural module interchanges these classes.

\item[{\em (f)}] $n=7$, $\pi\cap\pi(G)=\{2,3,5,7\}$,
$\vert G\vert_\pi=2^9\cdot 3^4\cdot5\cdot7$, and $H\simeq\Omega_7(2)$. There exist precisely two classes of conjugate subgroups of
this type in $G$, and  ${\rm SO}_7(q)$ interchanges these classes.

\item[{\em (g)}] $n=8$, $\eta=+$, $\pi\cap\pi(G)=\{2,3,5,7\}$, $\vert G\vert_\pi=2^{13}\cdot3^5\cdot5^2\cdot7$, and
$H\simeq 2\arbitraryext\Omega_8^+(2)$. There exist precisely four classes of conjugate subgroups of
this type in $G$. The subgroup of $\Out(G)$ generated by diagonal and graph automorphisms is isomorphic to $\Sym_4$ and acts on the set of
these classes as $\Sym_4$ in its natural permutation representation, and every diagonal automorphism acts without fixed points.

\item[{\em (h)}]\ $n=9$, $\pi\cap\pi(G)=\{2,3,5,7\}$,
$\vert G\vert_\pi=2^{14}\cdot3^5\cdot5^2\cdot7$, and $H\simeq 2\arbitraryext\Omega_8^+(2)\arbitraryext 2$. There exist precisely two classes
of conjugate subgroups of
this type in $G$, and  ${\rm SO}_9(q)$ interchanges these classes.
\end{itemize}
\item[{\em (B)}] Conversely, if one of the statements {\rm (a)--(h)} holds, then $G$ possesses a $\pi$-Hall subgroup with the given
structure.
\end{itemize}
\end{Lemma}

\begin{Lemma} \label{HallSubgroupsOfExceptionalGroups}  {\rm \cite[Lemmas~7.1--7.6]{NumbCl}, \cite[Theorem~8.13]{Surv}}   Assume that
$$G\in\{E_6^\eta(q),E_7(q),E_8(q),F_4(q),G_2(q),{{}^3D}_4(q)\},$$ where $q$  is a power of a prime $p$. Let  $\varepsilon=\varepsilon(q)$. Let $\pi$ be a
set of
primes such that  $2,3\in\pi$, $p\not\in\pi$.
 Then $G$ contains a $\pi$-Hall subgroup $H$  if and only if one of the following statements hold:
 \begin{itemize}
\item[{\em (a)}] $G$ is a group in Table~{\rm \ref{t2}} and  the values for the untwisted Lie rank $l$ of $G$, $\delta$ and the structure of the Weyl group $W$ are given in the
Table~{\rm \ref{t2}}; if $G=E_6^\eta(q)$ then $\eta=\varepsilon$; $\pi(W)\subseteq \pi\cap \pi(G)\subseteq \pi(q-\varepsilon)$, $H$ is a $\pi$-Hall subgroup of a
group
$T\arbitraryext W$, where $T$ is a maximal torus of order $({q-\varepsilon})^l/\delta$. All $\pi$-Hall subgroups of this type are conjugate in~$G$;

\item[{\em (b)}]  $G={{}^3D}_4(q)$,  $\pi\cap \pi(G)\subseteq
\pi(q-\varepsilon)$ and $H$  is a  $\pi$-Hall subgroup in
 $T\arbitraryext W(G_2)$, where $T$ is a maximal torus  of order $(q-\varepsilon)(q^3-\varepsilon)$. All $\pi$-Hall subgroups of this type
are conjugate in~$G$;

\item[{\em (c)}] $G=E_6^{-\varepsilon}(q)$, $\pi\cap \pi(G)\subseteq \pi(q-\varepsilon)$ and $H$ is a  $\pi$-Hall subgroup in
$T\arbitraryext W(F_4)$, where $T$ is a maximal torus of order ${(q^2-1)^2(q-\varepsilon)}^2/(3,q+\varepsilon)$. All $\pi$-Hall subgroups of this type
are conjugate in~$G$;

\item[{\em (d)}] $G=G_2(q)$, $\pi\cap \pi(G)=\{2,3,7\}$, $(q^2-1)_{\mbox{}_{\{2,3,7\}}}=24$, $(q^4+q^2+1)_{\mbox{}_7}=7$, $H\simeq G_2(2)$, and all
$\pi$-Hall
subgroups of
this type are conjugate in~$G$.
\end{itemize}
\end{Lemma}

\begin{longtable}{|l|c|c|c|r|}\caption{Weyl groups of exceptional root systems}\label{t2}\\ \hline
$G$ & $l$ &$\delta$ &$W$ & $|W|$\\
\hline\hline
$E^\eta_6(q)$& $6$  & $(3, q-\eta)$ & $W(E_6)\simeq{\rm Sp}_4(3)$ & $2^7.3^4.5$  \\
$E_7(q)$     & $7$  & $2$           & $W(E_7)\simeq 2\times{\rm P}\Omega_7(2)$ & $2^{10}.3^4.5.7$\\
$E_8(q)$     & $8$  & $1$           & $W(E_8)\simeq2\arbitraryext {\rm P}\Omega_8^+(2)\arbitraryext 2$  & $2^{14}.3^5.5^2.7$\\
$F_4(q)$     & $4$  & $1$           &$W(F_4)$ &$ 2^7.3^2$\\
$G_2(q)$     &  $2$ & $1$           &$W(G_2)$ &$2^2.3$\\
\hline
\end{longtable}

\begin{Lemma}
\label{p_ip_pi} {\rm \cite[Theorem~3.1]{Gross4}} Let  $G$ be a group of Lie
type with base field $\F_q$ of some  characteristic $p$. Assume that  $\pi$
is such that  $2,p\in\pi$,  and $3\notin\pi$. Suppose  $G\in \E_\pi$ and
$H\in\Hall_\pi(G)$. Then $p=2$ and one of the following statements holds.
\begin{itemize}
\item[$(1)$] $\pi\cap \pi(G)\subseteq\pi(q-1)\cup\{2\}$, a  Sylow
    $2$-subgroup $P$ of $H$ is normal in  $H$ and  $H/P$ is Abelian.
\item[$(2)$] $p=2$,  $G\simeq {}^2B_2(2^{2n+1})$ and $\pi(G)\subseteq\pi$.
\end{itemize}
\end{Lemma}

\begin{Lemma}\label{SolvabilityD2,3} {\rm \cite[Lemma 4]{GRV}} {\it  Let $G$ be
a nonabelian simple group. Then $G\in\D_{\{2,3\}}$ if and only if $G$ is a
Suzuki group ${}^2B_2(q)$. In this case every $\pi$-subgroup of $G$ is $2$-group.}
\end{Lemma}

\begin{Lemma}
\label{p_notip_pi}{\rm \cite[Lem\-ma~5.1 and Theorem~5.2]{Hall3'}}
Let  $G$ be a group of Lie type over a field
of characteristic  $p$.  Assume that  $\pi$  is such that $3,p\notin\pi$ and   $2\in\pi$.  Suppose $G\in \E_\pi$ and $H\in\Hall_\pi(G)$. Then either $H$ possesses a normal abelian $2'$-Hall subgroup or $G\cong {}^2G_2(3^{2m+1})$ and $\pi\cap\pi(G)=\{2,7\}$.
\end{Lemma}

\begin{Lemma}
\label{2or3notin_pi}
Let  $G$ be a simple nonabelian group. Assume that  $\pi$  is such that $\pi(G)\not\subseteq\pi$, $|\pi\cap\pi(G)|>1$, $2\in\pi$ and $3\notin\pi$. Suppose $G\in \E_\pi$ and $H\in\Hall_\pi(G)$. Then  the following statements hold.
\begin{itemize}
  \item[$(1)$] $H$ is solvable.
  \item[$(2)$] If every solvable $\pi$-subgroup of $G$ is conjugate to a subgroup of $H$, then  $G$ is a group of Lie type over a field of characteristic $p\notin\pi$ and $G\in\D_\pi$.
\end{itemize}
\end{Lemma}

\noindent{\sc Proof}. Statement~(1) is proved in \cite[Lemma~10]{GRV}.

Prove~(2).
Consider all possibilities for $G$, according to the classification of finite simple groups (see \cite[Theorem~0.1.1]{AschLyoSmSol}).

C\,a\,s\,e~1: $G\simeq \Alt_n$, $n\geq 5$. This case is impossible by Lemma~\ref{HallSymmetric}.

C\,a\,s\,e~2: $G$ is either a sporadic group or a Tits group. By Lemma \ref{HallSpor23} it follows that  $G\simeq J_1$ and $\pi\cap\pi(G)=\{2,7\}$. Now, by the Burnside theorem \cite[Ch.~I,~2]{DH}, every $\pi$-subgroup of $G$ is solvable and is conjugate to a subgroup of $H$, that is $G\in\D_\pi$. It contradicts Lemma~\ref{DpiCrit}.

C\,a\,s\,e~3: $G$ is a group of Lie type over a field $\mathbb{F}_q$ of characteristic $p\in\pi$.  By Lemma \ref{p_ip_pi}, $p=2$, $H$ is solvable and a Sylow $2$-subgroup $P$ of $H$ is normal in $H$. Moreover, $\pi\cap\pi(G)\subseteq\pi(q-1)\cup\{2\}$. Condition $|\pi\cap\pi(G)|>1$ implies that $q>2$. It is known that $G$ has a subgroup which is a homomorphic image of $\SL_2(q)=\PSL_2(q)$. Since $\PSL_2(q)$ is simple, we assume that $\PSL_2(q)\le G$. Take $r\in\pi\cap\pi(q-1)$. Then $\PSL_2(q)$ contains a dihedral subgroup $U$ of order $2r$ and $U$ has no normal Sylow 2-subgroups. Hence, $U$ is not conjugate to any subgroup of~$H$.

C\,a\,s\,e~4: $G$ is a group of Lie type over a field of characteristic $p\notin\pi$. Lemma \ref{p_notip_pi} implies that either $H$ possesses a normal abelian $2'$-Hall subgroup or $G\cong {}^2G_2(3^{2m+1})$, $\pi\cap\pi(G)=\{2,7\}$. In the last case every $\pi$-subgroup of $G$ is solvable by the Burnside theorem \cite[Ch.~I,~2]{DH} and $G\in\D_\pi$. Suppose, $H$ possesses a normal abelian $2'$-Hall subgroup. It is sufficient to prove that every $\pi$-subgroup of $G$ is solvable. Suppose $U$ is a nonsolvable subgroup of $G$. Then $U/U_{\mathfrak{S}}\ne 1$ and every minimal subnormal subgroup of $U/U_{\mathfrak{S}}$ is nonabelian simple and is isomorphic to a Suzuki group ${}^2B_2(2^{2m+1})$ in view of condition $3\notin\pi$ and the Thompson-Glauberman theorem \cite[Chapter II, Corollary 7.3]{Glauberman}. Take in $U$ the full preimage $V$ of a Borel subgroup of ${}^2B_2(2^{2m+1})$. Then $V$ is solvable and is  conjugate to a subgroup of $H$. In particular, $V$ and a Borel subgroup $V/U_{\mathfrak{S}}$ of ${}^2B_
2(2^{2m+1})$ possesses
a normal $2'$-Hall subgroup, but this is not true.
\qed\medskip

\begin{Lemma}\label{SolvabilityHallSubgroupsInDpi}  {\it Let $G\in\D_\pi$ be a nonabelian simple group. Then either $G$ is a $\pi$-group or every $\pi$-Hall subgroup of $G$ is solvable. In particular, if $G$ is not a $\pi$-group then $G\in\D_\tau$ for every $\tau\subseteq\pi$.}
\end{Lemma}

\noindent{\sc Proof}. Lemma follows from Lemmas~\ref{SolvabilityD2,3} and~\ref{2or3notin_pi} and the solvability of groups of odd order~\cite{FeitThompsonOddOrder}.
\qed\medskip

\subsection{Degrees of minimal faithful permutation representation}

In the following Lemma we collect some statements about minimal degrees of faithful permutation representations of some groups.

\begin{Lemma}\label{mu} The following statements hold.
\begin{itemize}
  \item[$(1)$] If $H\le G$ then $\mu(H)\le \mu(G)$.
  \item[$(2)$] {\rm \cite[Theorem 2]{Holt}} Let $G$ be a finite group. Let
      $\mathfrak{L}$ be a  complete class of finite groups. Let $N$ be the
      $\mathfrak{L}$-radical of $G$, that is the maximal normal
      $\mathfrak{L}$-subgroup of $G$. Then $\mu(G)\ge\mu(G/N)$.
  \item[$(3)$] {\rm \cite[Theorem 3.1]{Praeger}} If $G=L_1\times L_2\times\dots\times L_r$ and $L_1, L_2,\dots, L_r$ are simple then $$\mu(G)=\mu(L_1)+\mu(L_2)+\dots+\mu(L_r).$$
  \item[$(4)$] If $G$ is simple then $\mu(G)$ is equal to the minimum of indices of maximal subgroups in~$G$.
  \item[$(5)$] $\mu(\Sym_n)=\mu(\Alt_n)=n$.
\end{itemize}

\end{Lemma}

\subsection{Some subgroups of quasisimple and almost simple groups}

\begin{Lemma}\label{2A5} {\rm \cite[Tables~8.1 and 8.2]{Bray}} Assume that $q^2\equiv 1\pmod 5$ and $q$ is a power of an odd prime. Then $\SL_2(q)$ contains a subgroup isomorphic to $\SL_2(5)$ and $\PSL_2(q)$ contains a subgroup isomorphic to $\PSL_2(5)\cong \Alt_5$.
\end{Lemma}

\begin{Lemma}\label{2A6} {\rm \cite[Tables~8.8 and 8.10]{Bray}}  Assume that $q\equiv \eta\pmod 4$, where $q$ is a power of an odd prime and $\eta=\pm 1$. Then $\SL_4^\eta(q)$ contains a subgroup isomorphic to $4\circ2^{1+4}\arbitraryext \Alt_6$.
\end{Lemma}

\begin{Lemma}\label{IrrSubgr} {\rm \cite[Lemma~1.24]{DpiCl} } Let $l$ be an odd prime, $q>2$ be a power of a prime. Assume
$G = \langle D, x\rangle$, where $D$ is the group of all diagonal matrices of  determinant $1$, and
$$
x=\left(\begin{array}{ccccc}
    0 & 1 &  &  &  \\
     & 0 & 1 &  &  \\
     &  & \ddots &  &  \\
     &  &  & 0 & 1 \\
    1 &  &  &  & 0
  \end{array}\right)\in \SL_r(q).
$$
Then $|G|=(q-1)^{r-1}r$ and $G$ is absolutely irreducible.

\end{Lemma}

\begin{Lemma}\label{FundamentalSubgroups}{\it Let $q$ be a power of an odd prime. Let pair $(G^*,m)$, where $G^*$ is a quasisimple group and  $m$ is an positive integer, appear in the following list:
\begin{itemize}
  \item[$(1)$] $(G^*,m)=(\SL^\eta_n(q),[n/2])$,  $n>2$, $\eta=\pm$;
  \item[$(2)$] $(G^*,m)=(\Sp_{n}(q),n/2)$,  $n>2$ is even;
  \item[$(3)$] $(G^*,m)=(\Omega_{n}(q),2[n/4])$,  $n>5$ is odd;
  \item[$(4)$] $(G^*,m)=(\Omega^+_{n}(q),2[n/4])$,  $n>6$ is even;
  \item[$(5)$] $(G^*,m)=(\Omega^-_{n}(q),2[(n-1)/4])$,  $n>6$ is even;
  \item[$(6)$] $(G^*,m)=(E^\eta_{6}(q),4)$, $G^*$ is a quotient of the universal group by a central subgroup;
  \item[$(7)$] $(G^*,m)=(E_{7}(q),7)$, $G^*$ is a quotient of the universal group by a central subgroup;
  \item[$(8)$] $(G^*,m)=(E_{8}(q),8)$.
\end{itemize}
Then $G^*$ contains a collection $\Delta$ of  subgroups such that
\begin{itemize}
  \item[{\rm (a)}] every member of $\Delta$ is isomorphic to $\SL_2(q)$,
  \item[{\rm (b)}]  if $K^*,L^*\in\Delta$ are distinct then $[K^*,L^*]=1$, and
  \item[{\rm (c)}]  $|\Delta|=m$.
\end{itemize}}
\end{Lemma}

\noindent{\sc Proof.} Lemma follows from Aschbacher's theory of fundamental subgroups. Recall that, if $G^*$ is a group from the lemma, then $K^\ast$ is a {\it fundamental subgroup}, if $K^*$ is conjugate to a subgroup generated by a long root subgroup $U$ and its opposite $U^-$. Every fundamental subgroup of $G^*$ is isomorphic to $\SL_2(q)$.    Fix a Sylow 2-subgroup $S^*$ of $G^*$ and consider $\Delta={\rm Fun}(S^*)$ consisting of all fundamental subgroups $L^*$ such that $S^*\cap L^*$ is a Sylow 2-subgroup of $L^*$. It follows from \cite[6.2]{Asch1} that distinct elements of $\Delta$ elementwise commute and follows from  \cite[Theorem~2]{Asch2} that $|\Delta|=m$.
\qed\medskip

\begin{Lemma}\label{BigSubgroupsSn} {\it Let $G=\Sym_n$, where $n\ge 5$, and $G$ contains a subgroup $H$ having $\Alt_m$  as a homomorphic image for some $m\in\{n-1,n\}$. Then $H\cong \Alt_m$.}
\end{Lemma}

\noindent{\sc Proof.} Suppose, $m=n$. In this case $|G:H|\le 2$ and $H\in\{\Sym_n,\Alt_n\}$. Since $\Sym_n$ has no $\Alt_n$ as a homomorphic image, we have $H=\Alt_n$.

  Suppose, $m=n-1$. First of all, note the following well-known fact: every subgroup $H_0$ of $G$ of index $n$ is isomorphic to $\Sym_{n-1}$. Indeed, let $K$
be the kernel of the action of $G$ by right multiplication on the set $\Omega$ of right cosets of $H$ in~$G$. Then $K\le H_0$ and $|G:K|\ge |G:H_0|=n> 2$. Since $G$ has a unique minimal normal subgroup $\Alt_n$ and its index equals~2, we have $K=1$. Therefore, $G$ is embedded in $\Sym(\Omega)\cong \Sym_n$ and $G\cong \Sym(\Omega)$. Since $H_0$ is a point stabilizer in $G$, we have $H_0\cong \Sym_{n-1}$.

Now let $L$ be the kernel of an epimorphism $H\rightarrow \Alt_{n-1}$. We need to show that $L=1$. If not, then
$$|H|=|L||\Alt_{n-1}|\ge 2|\Alt_{n-1}|=(n-1)! \text{ and } |G:H|\le n.$$ Since $G$ has no proper subgroups of index less than $n$, except $\Alt_n$, we have either
$H=G=\Sym_n$, or $H\cong \Sym_{n-1}$. But $\Alt_{n-1}$ is not a homomorphic image of any of the groups $S_n$ or $\Sym_{n-1}.$   \qed\medskip

\begin{Lemma}\label{SenfNornSn} {\it Let $G=\Sym_n$. Then $H=N_G(H)$ if $H\leqslant G$ and $|G:H|$ is odd.}
\end{Lemma}

\noindent{\sc Proof.} Let $S$ be a Sylow $2$-subgroup of $H$. Then $S$ is a Sylow $2$-subgroup of $G$ and $S=N_G(S)$ by \cite[Lemma~4]{CF}. So, $S=N_G(S)\leq H$, and $H=N_G(H)$ by the Frattini Argument.\qed\medskip

\begin{Lemma}\label{HallsbgrpsInWreathProd}
{\it Let $H$ be a $\pi$-Hall subgroup of $G=L\wr \Sym_n$. Denote by $L_1\times \ldots\times L_n$ the base of the wreath product.
Assume that $L_i$ possesses a $\pi$-Hall subgroup that is isomorphic to $H\cap L_i$ and is not conjugate with $H\cap L_i$ in $L_i$.
Then $G$ possesses a $\pi$-Hall subgroup $K$ such that $H$ and $K$ have the same composition factors and are not conjugate in~$G$.}
\end{Lemma}

\noindent{\sc Proof}. We can assume for the simplicity that $i=1$. We set $A=L_1\times \ldots\times L_n$ and $H_1=H\cap L_1$ and denote by $K_1$ a subgroup of $L_1$  that is isomorphic to $H_1$
but is not conjugate to $H_1$. Note that $G$ acts on the set $\Omega=\{L_1,\dots, L_n\}$ via conjugation and $A$ is the kernel of this action. Moreover, it follows from the definition
of a wreath product that $N_G(L_1)=C_G(L_1)L_1$.

Renumberring $\{L_1,\ldots,L_n\}$, if necessary, we may choose a right
transversal $h_1=1,\dots, h_m$ of $N_H(L_1)$ in $H$ so that $L_1^{h_i}=L_i$. Then  $L^{h_i}\ne L^{h_j}$ if $i\ne j$. In particular, $m\le n$.  So
$\{L_1,\dots, L_m\}$ is an orbit of $H$ on $\Omega$. Thus both $\Delta= \{L_1,\dots, L_m\}$ and $\Gamma=\{L_{m+1},\dots, L_n\}$ are $H$-in\-va\-ri\-ant. Set
$$
K_i=\left\{
\begin{array}{rll}
                                     K_1^{h_i} & \text{ for } & i=1,\dots, m\quad\\
                                     H\cap L_i & \text{ for } & i=m+1,\dots, n
                                   \end{array}
                                   \right.
$$
and $K_0=\langle K_i\mid i=1,\dots,n\rangle=K_1\times\ldots\times K_n.$
By construction, $K_0\le A$ and $K_0\cong H\cap A\in\Hall_\pi(A)$.

We clame that for every $h\in H$ there exists $a\in A$ such that $K_0^h=K_0^a$.

Indeed, take $h\in H$. Then there exists $\sigma\in \Sym_n$ such that $L_i^h=L_{i\sigma}$ for $i=1,\ldots,n$. Since $\Delta$ and $\Gamma$ are both $H$-in\-va\-ri\-ant, we obtain that $i\sigma\in\{1,\ldots,m\}$ for $i=1,\ldots,m$ and $i\sigma \in\{m+1,\ldots,n\}$ for $i=m+1,\ldots,n$.

Take $i\le m$.
Then $h_ih=xh_{i\sigma}$ for some  $x\in N_H(L_1)$. In this case
$$K_i^h=K_1^{h_ih}=K_1^{xh_{i\sigma}}.$$ Since $x\in N_H(L_1)\le N_G(L_1)=C_G(L_1)L_1$ and $K_1\le L_1$,  $K_1^x=K_1^{b}$ for some $b\in L_1$. Set $b_i=b_1^{h_i}\in L_i$. Then we have
$$K_i^h=K_1^{xh_{i\sigma}}=K_1^{bh_{i\sigma}}=\left(K_1^{h_{i\sigma}}\right)^{b^{h_{i\sigma}}}=K_{i\sigma}^{b_{i\sigma}}.$$
Thus, we see that there are $b_1\in L_1,\dots, b_m\in L_m$ such that
$$K_i^h=K_{i\sigma}^{b_{i\sigma}}$$ for every $i\le m$.

Let $a=b_1\dots b_m$. We show that $K_0^h=K_0^a$. Indeed, we have seen that $K_i^h=K_{i\sigma}^{b_{i\sigma}}=K_{i\sigma}^{a}$ if $i\le m$.
If $i> m$ then  $$K_i^h= H\cap L_i^h= H\cap L_{i\sigma}=K_{i\sigma}=K_{i\sigma}^a,$$ since $a$ centralizes $K_{j}$ for all $j>m$. Hence,
$$
K_0^h=\langle K_i^h\mid i=1,\dots,n\rangle=\langle K_{i\sigma}^a\mid i=1,\dots,n\rangle=K_0^a.
$$

Now Lemma~\ref{Ext} implies that there is $K\in\Hall_\pi(HA)\subseteq  \Hall_\pi(G)$ such that $K_0=K\cap A$.

The groups $H$ and $K$ have the same composition factors, since $$K/ K\cap A\cong KA/A=HA/A\cong H/H\cap A\text{ and } K\cap A=K_0\cong H\cap A.$$

Suppose, $K=H^g$ for some $g\in G$.  Then the image of $g$ in $G/A$ normalizes $HA/A=KA/A$. Note that $2\in\pi$ in view of Lemma~\ref{Gross}. Therefore, the index of $HA/A$ in $G/A\cong \Sym_n$ is odd.
Lemma~\ref{SenfNornSn} implies that $g\in HA$, and so we may assume that $g\in A$. Thus $K\cap A=H^g\cap A$, i.e. $K_0=K\cap A=K_1\times \ldots\times K_n$ and $H\cap A=H_1\times \ldots\times H_n$ are conjugate in $A=L_1\times\ldots\times L_n$. In particular, $K_1$ and $H_1$ are conjugate in $L_1$,
a contradiction.
\qed

\section{
Proof of Theorem~\ref{DX_means_Dpi}}

Theorem~\ref{DX_means_Dpi} says that,  for  a finite simple group $G$, the following statements are equivalent:

\begin{itemize}
  \item[$(1)$] $G\in\D_\X$ and
  \item[$(2)$] either $G\in \X$ or $\pi(G)\nsubseteq\pi$ and $G\in\D_\pi$.
\end{itemize}

$(2)\Rightarrow (1)$. Obviously $\X\subseteq \D_\X$. So we  need to prove that if $\pi(G)\nsubseteq\pi$ and $G\in\D_\pi$ then $G\in\D_\X$. Lemma~\ref{SolvabilityHallSubgroupsInDpi} implies that every $\pi$-Hall subgroup (hence every $\pi$-subgroup  of $G$, since $G\in \D_\pi$) is solvable, thus it belongs to $\X$. On the other hand,  every $\X$-subgroup of $G$ is a $\pi$-subgroup and so is contained in a $\pi$-Hall subgroup (we again use $G\in\D_\pi$ here). Therefore, $\m_\X(G)=\Hall_\pi(G)$. Hence every two $\X$-ma\-xi\-mal subgroups of $G$ are conjugate, i.\,e. $G\in\D_\X$.

$(1)\Rightarrow (2)$. This implication is much harder to prove. Its proof requires case by case consideration and we organize it in a series of steps, and divide it in subsections.

\subsection{Proof of the implication $(1)\Rightarrow (2)$: general remarks}

Assume that $G\in\D_\X$ and $G\notin \X$. We need to show that $G\in \D_\pi$.

Lemma \ref{HallXSubgroups} implies that

\begin{itemize}
  \item[$(i)$]  {\it  $\m_\X(G)=\X\cap\Hall_\pi(G)=\Hall_\X(G)$. In particular, $G\in\E_\pi$ and all elements of $\Hall_\X(G)$ are conjugate.  }\end{itemize}

Suppose by contradiction that  $G\notin \D_\pi$. Then

\begin{itemize}
  \item[$(ii)$]  {\it There exists a $\pi$-subgroup of $G$ which does not belong to $\X$.} \end{itemize}

Otherwise the $\pi$-subgroups of $G$ are exactly the $\X$-subgroups, thus the $\pi$-maximal subgroups of $G$ are conjugate, i.e. $G\in \D_\pi$.

  \smallskip

  The inclusion $\mathfrak{S}_\pi\subseteq \X$ and $(ii)$ immediately imply

\begin{itemize}  \item[$(iii)$]  {\it  There exists a non-solvable $\pi$-subgroup in $G$. }\end{itemize}
The solvability of primary and biprimary groups  \cite[Ch.~I,~2]{DH} and $(iii)$ implies

\begin{itemize}
  \item[$(iv)$]  {\it $|\pi\cap\pi(G)|>2$.} \end{itemize}

  The Feit--Thompson theorem  \cite{FeitThompsonOddOrder} implies
\begin{itemize}
  \item[$(v)$]  {\it  $2\in\pi\cap\pi(G)$.}
\end{itemize}

Moreover, it follows from $(v)$ and Lemma~\ref{2or3notin_pi} that

\begin{itemize}  {\it
  \item[$(vi)$] $3\in\pi\cap\pi(G)$.}
\end{itemize}

Now we prove that

\begin{itemize}
  \item[$(vii)$]  {\it $G$ has no solvable $\pi$-Hall subgroups.}
\end{itemize}

Indeed, if $G$ has a solvable $\pi$-Hall subgroup $H$, then $H\in\X\cap\Hall_\pi(G)=\m_\X(G)$. In view of $(v)$, $(vi)$ and the Hall theorem, $H$ contains a $\{2,3\}$-Hall subgroup $H_0$ and  $H_0\in\Hall_{\{2,3\}}(G)$. Take an arbitrary $\{2,3\}$-subgroup $U$ in $G$. Since $U$ is solvable and in view of $(v)$ and $(vi)$ we have $U\in\X$. Now $G\in\D_\X$ implies that $U$ is conjugate to a subgroup of $H$. Moreover, the solvability of $H$ means that $U$ is conjugate to a subgroup of $H_0$ by the Hall theorem. Hence $G\in\D_{\{2,3\}}$, a contradiction with Lemma~\ref{SolvabilityD2,3}.

\smallskip

Now we exclude all possibilities for $G$, considering finite simple groups case by case, according to the Classification of the finite simple groups.

\subsection{Alternating groups}

The following statement follows from Lemma~\ref{An}.

\begin{itemize}
  \item[$(viii)$]  {\it $G$ is not isomorphic to an alternating group.}
\end{itemize}

\subsection{Sporadic groups and Tits group}

Now, exclude any possibilities for $G$ to be a sporadic group.

\begin{itemize}
  \item[$(ix)$]  {\it $G$ is not isomorphic to the Mathieu group $M_{11}$.}
\end{itemize}

Suppose, $G=M_{11}$. According to Lemma~\ref{HallSpor23} and Table~\ref{tb0} and in view of $(v)$--$(vii)$ it is sufficient
to consider the situation $\pi\cap\pi(G)=\{2,3,5\}$ and a Hall $\X$-subgroup $H$ of $G$ is $M_{10}=\Alt_6\nonsplitext2$. Take a $\{2,3\}$-Hall subgroup $U$ of $G$ (this group appears in Table~\ref{tb0}). Since $U$ is a solvable $\pi$-group, we have $U\in\X$. Now $G\in \D_\X$ implies that $U$ is conjugate to a subgroup of~$H$. But this means that $H$ and its unique nonabelian composition factor $\Alt_6$ satisfy~$\E_{\{2,3\}}$. A contradiction with Lemma~\ref{HallSymmetric}.

\begin{itemize}
  \item[$(x)$]  {\it $G$ is not isomorphic to the Mathieu group $M_{22}$.}
\end{itemize}

According to Lemma~\ref{HallSpor23} and Table~\ref{tb0}, if $G=M_{22}$ then an $\X$-Hall  subgroup $H$ of $G$ is isomorphic to $2^4:\Alt_6$. But $G$ contains \cite{Atlas} a maximal subgroup $U\cong 2^4:\Sym_5$ which is an $\X$-group and is not isomorphic to a subgroup of~$H$.

\begin{itemize}
  \item[$(xi)$]  {\it $G$ is not isomorphic to the Mathieu group $M_{23}$.}
\end{itemize}

Suppose, $G=M_{23}$ and $H\in\Hall_\X(G)$. Lemma~\ref{HallSpor23}, Table~\ref{tb0} and $(vii)$ imply that one of the following cases holds.
\begin{itemize}
  \item[(a)] $\pi\cap\pi(G)=\{2,3,5\}$ and $H\cong 2^4:\Alt_6$;
  \item[(b)] $\pi\cap\pi(G)=\{2,3,5\}$ and $H\cong 2^4:(3\times \Alt_5)$;
  \item[(c)] $\pi\cap\pi(G)=\{2,3,5,7\}$ and $H\cong \PSL_3(4):2_2$;
  \item[(d)] $\pi\cap\pi(G)=\{2,3,5,7\}$ and $H\cong 2^4:\Alt_7$;
   \item[(e)] $\pi\cap\pi(G)=\{2,3,5,7,11\}$ and $H\cong M_{22}$.
\end{itemize}

In Case (a), we consider an $\X$-subgroup $U\cong 2^4:(3\times \Alt_5)$ which is a $\pi$-Hall subgroup (and appears in Case (b)) and is not isomorphic to $H$.

\smallskip

Suppose, Case (b) holds. In $G$, consider a $\{2,3\}$-subgroup $U\cong 3^2:Q_8$, a Frobenius group which is contained in $\PSL_3(4)$, see~\cite{Atlas}. Suppose, $U$ is a subgroup of $H$. Let $$\overline{\phantom{x}}:H\rightarrow H/O_2(H)$$ be the natural epimorphism. Since $U$ has no non-trivial normal 2-subgroups, we have
$$
U\cong \overline{U}\le \overline{H}\cong 3\times \Alt_5.
$$
Now $|\overline U|_3=|\overline H|_3=3^2$, i.e. $\overline U$ contains a Sylow 3-subgroup of $\overline H$ and the cyclic subgroup $O_3(\overline H)$ of order 3 must be a normal subgroup in $\overline U$. But $U\cong \overline{U}$ has no normal subgroups of order~3. A contradiction.

\smallskip

We exclude Cases (c), (d) and (e), since the subgroup $H$ does not contain elements of order 15 in these cases while $M_{23}$ has a cyclic subgroup $U$ of order 15 and $U\in \mathfrak{S}_\pi\subseteq\X$.

\begin{itemize}
  \item[$(xii)$]  {\it $G$ is not isomorphic to the Mathieu group $M_{24}$.}
\end{itemize}

If $G=M_{24}$ then an $\X$-Hall subgroup $H$ is isomorphic to $2^6:3\nonsplitext \Sym_6$. Consider an $\X$ subgroup $U\cong 2^4:\Alt_6$ which is included in a maximal subgroup $M=M_{23}$ of $G$ and is a $\{2,3,5\}$-Hall subgroup of $M$. Since $G\in\D_\X$, without loss of generality, we can assume that $U\le H$. Now $U$ contains a subgroup $U_0\cong \Alt_6$ and, clearly, $U_0\cap O_2(H)=1$.  Let $$\overline{\phantom{x}}:H\rightarrow H/O_2(H)$$ be the natural epimorphism. We have
$$
\Alt_6\cong U_0\cong \overline{U}_0\le \overline{H}\cong 3\nonsplitext \Sym_6.
$$
But this means that $$\Alt_6\cong \overline{U}_0= \overline{U}_0'\le \overline{H}'\cong 3\nonsplitext \Alt_6$$ and we have a contradiction.

\begin{itemize}
  \item[$(xiii)$]  {\it $G$ is not isomorphic to the Janko group $J_{1}$.}
\end{itemize}

Suppose, $G=J_1$ and $H\in\Hall_\X(G)$. It follows from Lemma~\ref{HallSpor23}, Table~\ref{tb0} and $(vii)$ that $H\cong 2\times \Alt_5$. Clearly, $H$ contains no elements of order~15. But $G$ has a cyclic subgroup $U$ of order 15 (see \cite{Atlas}) and $U\in\mathfrak{S}_\pi\subseteq\X$.

\begin{itemize}
  \item[$(xiv)$]  {\it $G$ is not isomorphic to the Janko group $J_{4}$.}
\end{itemize}

Suppose, $G=J_1$ and $H\in\Hall_\X(G)$.  Lemma~\ref{HallSpor23} and Table~\ref{tb0}  imply that $$H\cong 2^{11}:2^6:3\nonsplitext \Sym_6.$$ We exclude this possibility arguing exactly as in $(xii)$, because $G$ contains a subgroup isomorphic to $M_{24}$.

\begin{itemize}
  \item[$(xv)$]  {\it $G$ is not isomorphic to any sporadic group or a Tits group.}
\end{itemize}

This statement follows from $(v)$, Lemma~\ref{HallSpor23}, and $(ix)$--$(xiv)$.

\subsection{Groups of Lie type of characteristic in  $p\in\pi(\X)$}

Now, according to Lemma~\ref{3.3.12}, we exclude the possibilities for $G$ to be isomorphic to a group of Lie type whose characteristic belongs to~$\pi$.

\begin{itemize}
  \item[$(xvi)$]  {\it If $G$ is a group of Lie type, then $G$ has no $\pi$-Hall subgroups contained in a Borel subgroup.}
\end{itemize}

Since every Borel subgroup of $G$ is solvable, $(xvi)$ follows from $(vii)$.

\begin{itemize}
  \item[$(xvii)$]  {\it $G$ is not isomorphic to $D_l(q)$, where $q$ is a power of some $p\in\pi$.}
\end{itemize}

Suppose, $G\cong D_l(q)$ and the numeration of the roots in a fundamental root system $\Pi$ of $G$ is chosen as in the Dynkin diagram on Pic.~$1$. It follows from Lemma~\ref{3.3.12} that $q$ is a power of~$2$,  $l$ is a Fermat prime (in particular, $l\ge 5$), and $(l,q-1)=1$. Moreover, if $H\in\Hall_\X(G)$, then $H$ is conjugate to the canonic parabolic maximal subgroup corresponding to the set $\Pi\setminus\{r_1\}$ of fundamental roots. This parabolic subgroup has a composition factor isomorphic to $D_{l-1}(q)$. Since  $\X$ is a complete class, we obtain that
$$D_i(q)\in \X \quad \text{ for } \quad i\le l-1, \quad \text{ and } \quad A_1(q)\in \X.$$
Moreover, $\pi(q-1)\subseteq\pi$. Consider the canonic parabolic maximal subgroup $P_J$ of $G$, corresponding to the set $J=\Pi\setminus\{r_2\}$. Above remarks and the completeness of $\X$ under extensions implies that $P_J\in\X$: the nonabelian composition factors of $P$ are isomorphic to $D_{l-2}(q)$ and, possibly, $A_1(q)$, while the orders of abelian composition factors belong to $\pi(q-1)\cup\{2\}\subseteq\pi$. But the maximality of $P$ means that $P_J$ is not conjugate to any subgroup of $H$, a contradictions with~$G\in \D_\X$.

\begin{itemize}
  \item[$(xviii)$]  {\it $G$ is not isomorphic to ${}^2D_l(q)$, where $q$ is a power of some $p\in\pi$.}
\end{itemize}

 Suppose, $G\cong {}^2D_l(q)$ and the numeration of the roots in a fundamental root system $\Pi^1$ of $G$ is chosen as in the Dynkin diagram on Pic.~$2$. It follows from Lemma~\ref{3.3.12} that $q$ is a power of~$2$,  $l-1$ is a Mersinne prime,  and $(l-1,q-1)=1$. Take $H\in\Hall_\X(G)$. Then $H$ is conjugate to the canonic parabolic maximal subgroup corresponding to the set $\Pi^1\setminus\{r^1_1\}$ of fundamental roots. This parabolic subgroup has a composition factor isomorphic to ${}^2D_{l-1}(q)$ if $l>4$ or isomorphic to ${}^2A_3(q)$ if $l=4$.  Consider the canonic parabolic maximal subgroup $P_J$ of $G$ which corresponds to the set $J=\Pi^1\setminus\{r_2^1\}$ of fundamental roots. Arguing as in $(xvii)$, we see that  $P_J\in\X$ and $P$ is not conjugate to any subgroup of $H$, a contradiction with $G\in \D_\X$.

\begin{itemize}
  \item[$(xix)$]  {\it $G$ is not isomorphic to $A_{l-1}(q)\cong \PSL_{l}(q)$, where $q$ is a power of some $p\in\pi$.}
\end{itemize}

Suppose, $G=\PSL_n(q)$, where $q$ is a power of some $p\in\pi$, and let $G^*=\SL_n(q)$. Lemma~\ref{Bijection} implies that $G\in\D_\X$ if and only if $G^*\in\D_\X$. Thus, $G^*\in\D_\X$ and, moreover, it follows from $(vii)$ and Lemma~\ref{Bijection} that there are no solvable $\pi$-Hall subgroups in $G^*$.

Identify $G^*$ with $\SL(V)$, where $V=\mathbb{F}_q^n$ is the natural $n$-dimensional module for $G^*$.  Let $H^*\in\Hall_\X(G^*)$.
By Lemma~\ref{3.3.12}, $H^*$ is the stabilizer in $G^*$ of a series
$$0={ V}_0<{ V}_1<\dots <{ V}_s={ V}$$ of subspaces such that ${\rm
dim}{ V}_i/{ V}_{i-1}=n_i$, $i=1,2,\dots ,s$, and one of the following conditions holds:
\begin{itemize}
\item[{\rm (a)}]  $n$~is a prime, $s=2$,
$n_1,n_2
\in
\{1,n-1\}$;

\item[{\rm (b)}] $n=4$,  $s=2$,
$n_1=n_2=2$; moreover, $q=2^{2t+1}$;

\item[{\rm (c)}] $n=5$, $s=2$,
$n_1,n_2\in\{2,3\}$;

\item[{\rm (d)}] $n=5$,  $s=3$,
$n_1,n_2,n_3\in\{1,2\}$;

\item[{\rm (e)}] $n=7$,  $s=2$,
$n_1,n_2\in\{3,4\}$;

\item[{\rm (f)}] $n=8$,  $s=2$,
$n_1=n_2=4$; moreover, $q=2^{2t}$;

\item[{\rm (g)}] $n=11$,  $s=2$, $n_1,n_2\in\{5,6\}$.
\end{itemize}

In cases (a), (c), (e),  and (g), $H^*$ is the stabilizer of a subspace of some dimension $m\ne n-m$ and the stabilizer $K^*$ of a subspace of dimension $n-m$ is isomorphic to $H^*$ (in particular, $K^*\in\X$) but is not conjugate to $H^*$. It contradicts $G^*\in \D_\X$.

 If case (d) holds, then there are exactly three conjugacy classes of $\pi$-Hall subgroups with the same composition factors and $H^*$ belongs to one of them. Thus, case (d) is impossible for  $G^*\in\D_\X$.

 Now consider cases (b) and (f). In these cases $n=4$ and $n=8$, respectively. Moreover, if $q=2$ then case (b) holds and $G=\PSL_4(2)\cong \Alt_8\notin \D_\X$ in view of $(viii)$. Therefore, we assume that $q>2$ if $n=4$.  Define  $r=n-1=3$ in (b) and  $r=n-1=7$ in (f). It is easy to check that  $r\in \pi$ in both cases. Consider the subgroup $U^*$ of $G^*$, consisting of all matrices of type
 $$
 \left(
 \begin{array}{cc}
   a &  \\
    & 1
 \end{array}
 \right),
 $$
where $a\in\langle D,x\rangle\le\SL_r(q)$, $D$ is the group of all diagonal matrices in $\SL_r(q)$ and
$$
x=\left(\begin{array}{ccccc}
    0 & 1 &  &  &  \\
     & 0 & 1 &  &  \\
     &  & \ddots &  &  \\
     &  &  & 0 & 1 \\
    1 &  &  &  & 0
  \end{array}\right)\in \SL_r(q).
$$
By Lemma~\ref{IrrSubgr} it follows that there is a subspace $W$ of $V$ of dimension $r$, such that $U^*$ acts irreducibly on~$W$. Clearly, $U^*$ cannot stabilize any subspace of dimension $n/2=2$ in case (b) or $n/2=4$  in case (f). Therefore, $U^*$ is not conjugate to any subgroup of~$H^*$. By Lemma~\ref{IrrSubgr}, $U^*$ is a solvable $\pi$-group, so $U^*\in \X$, a contradiction with  $G^*\in\D_\X$.

\begin{itemize}
  \item[$(xx)$]  {\it $G$ is not isomorphic to any group of Lie type of characteristic $p\in\pi$.}
\end{itemize}

This statement follows from $(xvi)$--$(xix)$ and Lemma~\ref{3.3.12}.

\subsection{Classical groups of characteristic  $p\notin\pi(\X)$ }

In view of $(xx)$, $G$ is a group  of Lie type over a field of an order $q$ and  characteristic $p\notin\pi$. In particular, $p\ne 2,3$.

We start with the smallest case $G=\PSL_2(q)$.

\begin{itemize}
  \item[$(xxi)$]  {\it $G$ is not isomorphic to $\PSL_2(q)$.}
\end{itemize}
Suppose $G=\PSL_2(q)$, and denote $G^*=\SL_2(q)$. Then $G^*\in \D_\X$ and  $G^*$ has no solvable $\pi$-Hall subgroups  by $(vii)$ and Lemmas~\ref{HallSubgroup} and~\ref{Bijection}, so statement (d) of Lemma~\ref{SLU23dim2} holds. Therefore if $H^*$ is an $\X$-Hall subgroup of $G^*$ then the image of $H^*$ in $G^*/Z(G^*)\cong G$ is isomorphic to $\Alt_5$. But in this case there are exactly two conjugacy classes of $\X$-Hall subgroup in $G$. It contradicts $G\in\D_\X$.

\medskip

Now we show that  $G$ is not isomorphic to a classical group. First we consider the most transparent case of symplectic groups.
Similar, but more complicated, arguments appear in the consideration of the other types of classical groups: linear, unitary and orthogonal.

 \begin{itemize}
  \item[$(xxii)$]  {\it $G$ is not isomorphic to $\PSp_{2n}(q)$.}
\end{itemize}

Suppose $G=\PSp_{2n}(q)$ and denote $G^*=\Sp_{2n}(q)$. By $(vii)$ and Lemma~\ref{Bijection}, we have $G^*\in \D_\X$ and $G^*$ has no solvable $\pi$-Hall subgroups.
Consider $H^*\in\Hall_\X(G^*)$. We claim  that
\begin{itemize}
    \item  $\pi\cap\pi(G^*)\subseteq\pi(q^2-1)$;
    \item $H^*$ is included in a subgroup $$M^*\cong \SL_2(q)\wr \Sym_n,$$ we denote by $B^*$   the base of this wreath product;
    \item $H^*/(H^*\cap B^*)$ is isomorphic to a $\pi$-Hall subgroup of $\Sym_n$;
    \item $H^*\cap B^*$ is solvable.
\end{itemize}

First two items can be found in Lemma \ref{HallSubgroupsOfSymplecticGroups}, the third item follows by Lemma \ref{HallSubgroup}.
The last item follows by Lemma \ref{HallsbgrpsInWreathProd}, since if $H^*\cap B^*$ is nonsolvable, then, for some component $L^*=SL_2(q)$, $H^*\cap L^*$ is a
nonsolvable $\pi$-Hall subgroup of $L^*$. Now Lemma \ref{SLU23dim2} implies that $L^*$ possesses $\pi$-Hall subgroup that is isomorphic and nonconjugate to $H^\ast\cap  L^*$.
Finally Lemma \ref{HallsbgrpsInWreathProd} implies that $M^*$ possesses a $\pi$-Hall subgroup that  is nonconjugate to $H^*$ but have the same composition factors.
Lemma \ref{HallSubgroupsOfSymplecticGroups}(C) implies that $G^*$ possesses nonconjugate $\X$-Hall subgroups, a contradiction with $G^*\in \D_\X$.

The nonsolvability  of $H^*$ and Lemma~\ref{HallSymmetric} imply that $H^*/(H^*\cap B^*)$ is isomorphic to a symmetric group of degree $n$ or $n-1$ and this degree is  at least $5$. In particular, $5\in\pi\cap\pi(G^*)$, $\Alt_5\in \X$, and $5$ divides $q^2-1$. Moreover, $H^*\cap B^*$ coincides with the solvable radical $H^*_{\mathfrak{S}}$ of~$H^*$.

Lemma~\ref{FundamentalSubgroups} implies that  $G^*$ possesses a collection $\Delta$ of subgroups isomorphic to $\SL_2(q)$ such that $|\Delta|=n$ and $[K^*, L^*]=1$ for every $K^*,L^*\in\Delta$  and $K^\ast\not=L^\ast$. Since $5$ divides $q^2-1$,  Lemma~\ref{2A5} implies that  $\SL_2(q)$ possesses a subgroup isomorphic to $\SL_2(5)$. For every $K^\ast\in \Delta$ fix some  $U(K^*)\le K^*$ such that $U(K^*)\cong \SL_2(5)$. Set
 $$
 U^*=\langle U(K^*)\mid K^*\in \Delta \rangle.
 $$
 It follows from the definition that $$U^*/U^*_{\mathfrak{S}}\cong \underbrace{\Alt_5\times\dots\times \Alt_5}_{n \text{ times}}$$
 and $U^*_{\mathfrak{S}}$ is a 2-group. Thus, $U^*\in\X$.

We show that $U^*$ is not conjugate to a subgroup of $H^*$, and this  contradicts $G^*\in\D_\X$. Indeed, if $U^\ast$ is conjugate to a subgroup of $H^\ast$, then we can
assume that $U^*\le H^*$. Denote by $R^*=H^*\cap B^*$ the solvable radical of $H^*$ and let $$\overline{\phantom{x}}:H^*\rightarrow H^*/R^*=\overline{H}{}^*$$ be the natural epimorphism.
We have seen above that $\overline{H}{}^*$ is isomorphic to a subgroup of  $\Sym_n$. Therefore, $$\mu(\overline{U}{}^*)\le \mu(\overline{H}{}^\ast)\le n.$$
On the other hand, $\overline{U}{}^*/\overline{U}{}^*_{\mathfrak{S}}\cong U^*/U^*_{\mathfrak{S}}$ and Lemma~\ref{mu} implies that
$$\mu(\overline{U}{}^*)\ge \mu(U^*/U^*_{\mathfrak{S}})=5n>n.$$ It contradicts the previous inequality.

Thus, $(xxii)$ is proved.

\begin{itemize}
  \item[$(xxiii)$]  {\it $G$ is not isomorphic to $\PSL^\eta_n(q)$, $\eta=\pm$.}
\end{itemize}

Suppose $G=\PSL^\eta_n(q)$ and denote $G^*=\SL^\eta_n(q)$. By $(vii)$ and Lemma~\ref{Bijection}, we have $G^*\in \D_\X$ and $G^*$ has no solvable $\pi$-Hall subgroups. Let $H^*\in\Hall_\X(G^*)$. Consider all possibilities for $H^*$ given in statements (a)--(e) of Lemma~\ref{HallSubgroupsOfLinearAndUnitaryGroups}.

In  case (a)  $n$ is equal to $2$, and this case is excluded in view of $(xxi)$.

In  case (d)  $H^*$ is isomorphic to $4\arbitraryext  2^4\arbitraryext  \Alt_6$. In this case $G^*$ has two conjugacy classes of $\pi$-Hall subgroups isomorphic $4\arbitraryext  2^4\arbitraryext  \Alt_6$. So if $H^*$ satisfies (d), then this  contradicts~${G^*\in\D_\X}$.

In case (e) we have $\pi\cap\pi(G^*)=\{2,3\}$. So $H^\ast$ is solvable and this case is excluded  in view of $(vii)$.

Thus, one of the following statements holds.

\begin{itemize}
  \item[{ (b)}]
 $q\equiv \eta \pmod 4$, $\Sym_n$ satisfies $\E_\pi$, $\pi\cap\pi(G)\subseteq \pi(q-\eta)\cup \pi(n!)$, and if $r\in(\pi\cap\pi(n!))\setminus
\pi(q-\eta)$, then $|G^*|_r=|\Sym_n|_r$. In this case $H^*$ is included in  $$M^*=L^*\cap G^*\simeq (q-\eta)^{n-1}\arbitraryext \Sym_n,$$ where $L^*=\GL_1^\eta(q)\wr
\Sym_n\leq \GL_n^\eta(q)$.

\item[{ (c)}] $n=2m+k$, where $k\in\{0,1\}$, $m\ge 1$,
$q\equiv -\eta \pmod {3}$, $\pi\cap\pi(G)\subseteq \pi(q^2-1)$, the groups $\Sym_m$ and $\GL_2^\eta(q)$ satisfy $\E_\pi$.
In this case $H^*$ is contained in $$M^*=L^*\cap
G^*\simeq(\underbrace{\GL_2^\eta(q)\circ\dots\circ\GL_2^\eta(q)}_{\displaystyle
m \mbox{ \rm times}})\arbitraryext\Sym_m\circ Z,$$ where $L^*=\GL_2^\eta(q)\wr \Sym_m\times Z\leq \GL_n(q)$ and $Z$ is a cyclic group of order $q-\eta$ for
$k=1$, and $Z$ is trivial for
$k=0$. The intersection of $H^*$ with each factor $\GL_2^\eta(q)$ is a  $\pi$-Hall subgroup in $\GL_2^\eta(q)$.
\end{itemize}

By Lemma~\ref{GLU23dim2} $\pi$-Hall subgroups of $\GL_2^\eta(q)$ are solvable. Since $H^*$ is nonsolvable, every nonabelian composition factor of $H^*$ is a composition factor of a $\pi$-Hall subgroup of a symmetric group of degree at most $n$. In both cases (b) and (c), it follows from Lemma~\ref{HallSymmetric} that
 \begin{itemize}
   \item every nonabelian composition factor of $H^*$ is isomorphic to an alternating group; in particular
   \item $5\in\pi\cap\pi(G^*)$, $\Alt_5\in\X$ and $n\ge 5$; and
   \item $H^*/H^*_\mathfrak{S}$ is isomorphic to a subgroup of $\Sym_n$.
 \end{itemize}
Now we consider two cases:  $5$ divides $q^2-1$ and
 $5$ does not divide $q^2-1$.

Suppose,  $5$ divides $q^2-1$. In this case we argue similarly the case of symplectic groups above. Lemma~\ref{FundamentalSubgroups} implies that  $G^*$ possesses a collection $\Delta$ of subgroups isomorphic to $\SL_2(q)$ such that $|\Delta|=[n/2]$ and $[K^*, L^*]=1$ for every $K^*,L^*\in\Delta$  and $K^\ast\not=L^\ast$. Since $5$ divides $q^2-1$,  Lemma~\ref{2A5} implies that  $\SL_2(q)$ possesses a subgroup isomorphic to $\SL_2(5)$. For every $K^\ast\in \Delta$ fix some  $U(K^*)\le K^*$ such that $U(K^*)\cong \SL_2(5)$. Set
 $$
 U^*=\langle U(K^*)\mid K^*\in \Delta \rangle.
 $$
 It follows from the definition that $$U^*/U^*_{\mathfrak{S}}\cong \underbrace{\Alt_5\times\dots\times \Alt_5}_{[n/2] \text{ times}}$$
 and $U^*_{\mathfrak{S}}$ is a 2-group. Thus, $U^*\in\X$.

We claim that $U^*$ is not conjugate to a subgroup of $H^*$ and this  contradicts $G^*\in\D_\X$. Indeed, if $U^\ast$ is conjugate to a subgroup of $H^\ast$, then we can
assume that $U^*\le H^*$. Denote by $R^*$ the solvable radical of $H^*$ and let $$\overline{\phantom{x}}:H^*\rightarrow H^*/R^*=\overline{H}{}^\ast$$ be the natural epimorphism.
It follows from above and from Lemma~\ref{HallSymmetric} that $\overline{H}{}^*$ is isomorphic to a symmetric group of degree at most $n$.
Therefore, $$\mu(\overline{U}{}^*)\le n.$$ On the other hand, Lemma~\ref{mu} implies that  $$\mu(\overline{U}{}^*)\ge \mu(U^*/U^*_{\mathfrak{S}})=5[n/2]>n.$$ It contradicts the previous inequality. Hence  $5$ does not divide $q^2-1$.

 Suppose,  $5$ does not divide $q^2-1$. It means that case (b) holds (in particular, the solvable radical $H^*_{\mathfrak{S}}$ of $H^*$ is abelian) and $|G^*|_5=|\Sym_n|_5$. We have
 $$
 |G^*|_5=\prod\limits_{i=1}^n(q^i-\eta^i)_5 \quad \text{and}\quad |\Sym_n|_5=(n!)_5.
 $$
 Lemma~\ref{arithm} implies that $[n/4]=[n/5]$. Since $n\ge 5$, this means $$n\in\{5,6,7,10,11,15\}.$$

Assume that $n\in\{5,6,7\}$ first. Since $\Sym_n\in\E_\pi$ and a $\pi$-Hall subgroup of $\Sym_n$ belongs to $\X$, it follows from Lemma~\ref{HallSymmetric} that $\Alt_5\in\X$. Moreover, if $n=6$ or $n=7$, then $\Alt_6\in \X$.

The group $G^*$ has a subgroup isomorphic to $\SL_4^\eta(q)$. Moreover, (b) implies that $q\equiv \eta\pmod 4$ and by Lemma~\ref{2A6},  $G^*$ has a subgroup $$W^*\cong 4\circ 2^{1+4}. \Alt_6.$$
Define $U^\ast\leq G^\ast $ in the following way. If $n=6,7$, then $U^*=W^*.$ If $n=5$, then $W^*/W_{\mathfrak{S}}^*\cong \Alt_6$ contains a subgroup isomorphic to $\Alt_5$, and we set $U^*$ to be equal to its full preimage in $W^\ast$. By construction $U^*\in \X$.

We claim that $U^*$ is not conjugate to any subgroup of $H^*$ and this  contradicts $G^*\in\D_\X$. Indeed, if $U^*\le H^*$ and $R^*=H^*_{\mathfrak{S}}$ is the solvable radical of $H^*$ then $U^*/(U^*\cap R^*)$ is isomorphic to a subgroup of $H^*/R^*\leq \Sym_n$. We have that $U^*/U_{\mathfrak{S}}^*\cong \Alt_m$ for some $m\in \{n,n-1\}$. Since $U^*/U_{\mathfrak{S}}^*$ is a homomorphic image of $U^*/(U^*\cap R^*)$, it follows by Lemma~\ref{BigSubgroupsSn} that $U^*/(U^*\cap R^*)\cong \Alt_m$. Therefore, $U_{\mathfrak{S}}^*=U^*\cap R^*$. This is impossible, since $R^*$ is abelian, while $U_{\mathfrak{S}}^*\cong 4\circ 2^{1+4}$ contains an
extra special 2-subgroup of order~$2^5$.

Assume finally that $n\in\{10,11,15\}$. Lemma~\ref{HallSymmetric} implies that ${\Alt_{10}\in \X}$. Therefore $\Alt_6\in\X$. It is clear that $G^*$ has a subgroup
 $$ \SL_4^\eta(q)\circ \SL_4^\eta(q) \text{ if } n\in\{10,11\}, \text{ and } \SL_4^\eta(q)\circ\SL_4^\eta(q)\circ\SL_4^\eta(q) \text{ if } n=15.
 $$

By Lemma~\ref{2A6} and in view of $q\equiv \eta\pmod 4$, we can find a subgroup $U^*$ in $G^*$ such that $U_\mathfrak{S}^*=O_2(U^*)$ and   $$U^*/U_\mathfrak{S}^*\cong
\left\{
\begin{array}{rll}
                                     \Alt_6\times \Alt_6, & \text{ if } & n\in\{10,11\}, \quad\\
                                     \Alt_6\times \Alt_6\times \Alt_6, & \text{ if } & n=15.
                                   \end{array}
                                   \right.
$$
Clearly, $U^*\in\X$. But $U^*$ is not conjugate to a subgroup of $H^*$. Indeed, if $U^*\le H^*$, then $U^*/(U^*\cap R^*)$ is isomorphic to a subgroup of $\Sym_n$, where $R^*=H^*_{\mathfrak{S}}$. Therefore, by Lemma~\ref{mu} we have
$$
n\ge \mu(U^*/(U^*\cap R^*))\ge \mu(U^*/U^*_{\mathfrak{S}})=\left\{
\begin{array}{rll}
                                     \mu(\Alt_6\times \Alt_6)=12, & \text{ if } & n\in\{10,11\}, \quad \\
                                     \mu(\Alt_6\times \Alt_6\times \Alt_6)=18, & \text{ if } & n=15,
                                   \end{array}
                                   \right.
$$
a contradiction.

 Thus, $(xxiii)$ is proven.

\begin{itemize}
  \item[$(xxiv)$]  {\it $G$ is not isomorphic to $\P\Omega^\eta_n(q)$, $\eta\in\{+,-,\circ\}$.}
\end{itemize}

Suppose $G=\P\Omega^\eta_n(q)$, $n\ge 7$ and denote $G^*=\Omega^\eta_n(q)$. By $(vii)$ and Lemma~\ref{Bijection}, we have $G^*\in \D_\X$,  and $G^*$ has no solvable $\pi$-Hall subgroups and has exactly one class of $\X$-Hall subgroups. Let $H^*\in\Hall_\X(G^*)$. Consider all possibilities for $H^*$ given in statements (a)--(h) of Lemma~\ref{HallSubgroupsOfOrthogonalGroupsOfOddDimension}.

In cases (d) and (e) we have $\pi\cap\pi(G^*)=\{2,3\}$, and we exclude these cases in view of $(vii)$ and the solvability of $\{2,3\}$-groups.

We exclude cases (f), (g) and (h),   since in all these cases  there are at least two conjugacy classes of $\X$-Hall subgroups of $G^*$ isomorphic to $H^*$.

Thus, one of the following statements holds.

\begin{itemize}
\item [{\rm (a)}] $n=2m+1$, $\pi\cap\pi(G^*)\subseteq\pi(q-\varepsilon)$, $q\equiv \varepsilon\pmod {12}$, $\Sym_m\in \E_\pi$, and $H^*$ is a
$\pi$-Hall subgroup in $$M^*=\big(\Oo_2^\varepsilon(q)\wr \Sym_m\times \Oo_1(q)\big)\cap G^*.$$

\item [{\rm (b)}] $n=2m$, $\eta=\varepsilon^m$, $\pi\cap\pi(G^*)\subseteq\pi(q-\varepsilon)$, $q\equiv \varepsilon\pmod {12}$, $\Sym_m\in
\E_\pi$, and $H$ is a $\pi$-Hall
subgroup in $$M^*=\big(\Oo_2^\varepsilon(q)\wr \Sym_m\big)\cap G^*.$$

\item [{\rm (c)}] $n=2m$, $\eta=-\varepsilon^m$, $\pi\cap\pi(G^*)\subseteq\pi(q-\varepsilon)$, $q\equiv \varepsilon\pmod {12}$,
$\Sym_{m-1}\in \E_\pi$, and $H^*$ is a
$\pi$-Hall subgroup of $$M^*=\big(\Oo_2^\varepsilon(q)\wr \Sym_{m-1}\times \Oo^{-\varepsilon}_2(q)\big)\cap G^*.$$
\end{itemize}
Here $\varepsilon=\pm1$ and $q-\varepsilon$ is divisible by~4.

Groups $\Oo_2^+(q)$ and $\Oo_2^-(q)$ are solvable. As in the proofs of $(xxii)$ and $(xxiii)$, we see that the symmetric group of degree $m$, in cases (a) and (b), and of degree $m-1$ in case (c) has nonsolvable $\X$-Hall subgroup which is isomorphic to a symmetric group. Moreover, this $\X$-Hall subgroup is isomorphic to $H^*/H^*_{\mathfrak{S}}$. Thus,
\begin{itemize}
   \item $5\in\pi\cap\pi(G^*)\subseteq\pi(q-\varepsilon)\subseteq\pi(q^2-1)$,
   \item $\Alt_5\in\X$ and
   \item $H^*/H^*_\mathfrak{S}$ is isomorphic to a subgroup of $\Sym_m$  in cases (a) and (b) and of  $\Sym_{m-1}$  in case (c). Therefore, $$\mu(H^*/H^*_\mathfrak{S})\le m= [n/2].$$
 \end{itemize}

Lemma~\ref{FundamentalSubgroups} implies that  $G^*$ possesses a collection $\Delta$ of subgroups isomorphic to $\SL_2(q)$ such that $|\Delta|=k\ge 2[(n-1)/4]$ and $[K^*, L^*]=1$ for every $K^*,L^*\in\Delta$  and $K^\ast\not=L^\ast$. Since $5$ divides $q^2-1$,  Lemma~\ref{2A5} implies that  $\SL_2(q)$ possesses a subgroup isomorphic to $\SL_2(5)$. For every $K^\ast\in \Delta$ fix some  $U(K^*)\le K^*$ such that $U(K^*)\cong \SL_2(5)$. Set
 $$
 U^*=\langle U(K^*)\mid K^*\in \Delta \rangle.
 $$
 It follows from the definition that $$U^*/U^*_{\mathfrak{S}}\cong \underbrace{\Alt_5\times\dots\times \Alt_5}_{k \text{ times}}$$
 and $U^*_{\mathfrak{S}}$ is a 2-group. Thus, $U^*\in\X$.

We show that $U^*$ is not conjugate to a subgroup of $H^*$. Otherwise we can
assume that $U^*\le H^*$. Let $R^*=H^*_\mathfrak{S}$  and let $$\overline{\phantom{x}}:H^*\rightarrow H^*/R^*=\overline{H}{}^\ast$$ be the natural epimorphism.
Therefore, $$\mu(\overline{U}{}^*)\le\mu(\overline{H}{}^*)\le  [n/2].$$ On the other hand, since $n\ge 7$ and in view of Lemma~\ref{mu}, we have
$$\mu(\overline{U}{}^*)\ge \mu(U^*/U^*_{\mathfrak{S}})=5k\ge 10\left[\frac{n-1}{4}\right]\ge \frac{10(n-4)}{4}=\frac{5(n-4)}{2}>\frac{n+1}{2}\ge \left[\frac{n}{2}\right],$$
a contradiction.
\begin{itemize}
  \item[$(xxv)$]  {\it $G$ is not isomorphic to a classical group.}
\end{itemize}

This statement follows from $(xx)$ if characteristic of a group belongs to $\pi$ and from $(xxii)$--$(xxiv)$ in over cases.

\subsection{Exceptional groups of Lie type of characteristic  $p\notin\pi(\X)$ }

\begin{itemize}
  \item[$(xxv)$]  {\it $G$ is not isomorphic to one of groups ${}^2B_2(2^{2m+1})$,  ${}^2G_2(3^{2m+1})$, and  ${}^2F_4(2^{2m+1})$.}
\end{itemize}

This statement follows from $(v)$, $(vi)$ and  $(xx)$.

\begin{itemize}
  \item[$(xxvi)$]  {\it $G$ is not isomorphic to   $G_2(q)$.}
\end{itemize}

Suppose that $G=G_2(q)$ and $H\in\Hall_\X(G)$. By Lemma~\ref{HallSubgroupsOfExceptionalGroups}, either $H$ is solvable, which contradicts~$(vii)$, or statement (d) of Lemma~\ref{HallSubgroupsOfExceptionalGroups} holds:
\begin{itemize}
\item[{\rm (d)}] $G=G_2(q)$, $\pi\cap \pi(G)=\{2,3,7\}$, $(q^2-1)_{\mbox{}_{\{2,3,7\}}}=24$, $(q^4+q^2+1)_{\mbox{}_7}=7$, and $H\simeq G_2(2)$.
\end{itemize}
By~\cite{Atlas}, $H'\cong \PSU_3(3)$ has a maximal subgroup isomorphic to $\SL_3(2)\cong \PSL_2(7)$ and every maximal subgroup of $H'$ not isomorphic to $\SL_3(2)$ is solvable. This implies that $\SL_3(2)\in \X$ and $H$ has no subgroups isomorphic to $2^3\arbitraryext \SL_3(2)$, which belongs to $\X$. On the other hand, it follows from~\cite[Table~1]{CLSS} that $G$ has a subgroups isomorphic to $2^3\arbitraryext \SL_3(2)$.

\begin{itemize}
  \item[$(xxvii)$]  {\it $G$ is not isomorphic to one of groups   ${}^3D_4(q)$ and  $F_4(q)$.}
\end{itemize}

By Lemma~\ref{HallSubgroupsOfExceptionalGroups}, every Hall $\X$-subgroup of  ${}^3D_4(q)$ and  $F_4(q)$ is solvable, which contradicts~$(vii)$, if $G\in\{{}^3D_4(q),F_4(q)\}$.

\begin{itemize}
  \item[$(xxviii)$]  {\it $G$ is not isomorphic to one of groups     $E_6(q)$ and  ${}^2E_6(q)$.}
\end{itemize}

Suppose $G=E^\eta_6(q)$, $\eta=\pm$ and $H\in\Hall_\X(G)$. Since $H$ is not solvable, statement (c) of Lemma~\ref{HallSubgroupsOfExceptionalGroups} does not hold, and we have case (a) of this Lemma for $E^\eta_6(q)$:
\begin{itemize}
  \item  $4$ divides $q-\eta$, $\{2,3,5\}\subseteq \pi\cap \pi(G)\subseteq \pi(q-\eta)$, $H$ is a $\pi$-Hall subgroup of a
group
$T\arbitraryext \Sp_4(3)$, where $T$ is a maximal torus of order $({q-\eta})^6/3$.
\end{itemize}
Note that $\Sp_4(3)$ is a $\pi$-group. This implies that $\Sp_4(3)$ is a homomorphic image of $H$ and $\Sp_4(3)\in\X$. Furthermore, $H/H_{\mathfrak{S}}\cong \Sp_4(3)$. By information in~\cite{Atlas}, $\Sp_4(3)$ has a subgroup isomorphic to $\Alt_5$. Therefore, $\Alt_5\in\X$.

Lemma~\ref{FundamentalSubgroups} implies that  $G$ possesses a collection $\Delta$ of subgroups
isomorphic to $\SL_2(q)$ such that $|\Delta|=4$ and $[K, L]=1$ for every $K,L\in\Delta$  and $K\not=L$.
Since $5$ divides $q^2-1$,  Lemma~\ref{2A5} implies that  $\SL_2(q)$ possesses a subgroup isomorphic to
$\SL_2(5)$. For every $K\in \Delta$ fix some  $U(K)\le K$ such that $U(K)\cong \SL_2(5)$. Set
 $$
 U=\langle U(K)\mid K\in \Delta \rangle.
 $$
 It follows from the definition that $$U/U_{\mathfrak{S}}\cong \Alt_5\times\Alt_5\times\Alt_5\times \Alt_5$$
 and $U_{\mathfrak{S}}$ is a 2-group. Thus, $U\in\X$. Suppose that $U$ is conjugate to a subgroup of $H$. Then $H/H_{\mathfrak{S}}\cong \Sp_4(3)$  contains a subgroup for which $U/U_{\mathfrak{S}}$ is a homomorphic image. But
 $$|H/H_{\mathfrak{S}}|_5= |\Sp_4(3)|_5=5<5^4=|\Alt_5|_5^4=|U/U_{\mathfrak{S}}|_5,$$ and this is impossible.

 \begin{itemize}
  \item[$(xxix)$]  {\it $G$ is not isomorphic to      $E_7(q)$.}
\end{itemize}

Suppose $G=E_7(q)$ and $H\in\Hall_\X(G)$. By Lemma~\ref{HallSubgroupsOfExceptionalGroups} we have:
\begin{itemize}
  \item  $\{2,3,5,7\}\subseteq \pi\cap \pi(G)\subseteq \pi(q-\varepsilon)$, where $\varepsilon=\pm1$ is such that $4$ divides $q-\varepsilon$, $H$ is a $\pi$-Hall subgroup of a
group
$T\arbitraryext \big(2\times\P\Omega_7(2)\big)$, where $T$ is a maximal torus of order $({q-\eta})^7/2$.
\end{itemize}
This implies that $\P\Omega_7(2)\in\X$ and $H/H_{\mathfrak{S}}\cong \P\Omega_7(2)$. In $\P\Omega_7(2)$ there is a maximal subgroup $\Omega^+_6(2)\cong\Sym_8$. Therefore, $\Alt_5\in\X$.

Now we argue as in $(xxviii)$. Lemma~\ref{FundamentalSubgroups} implies that  $G$
possesses a collection $\Delta$ of subgroups isomorphic to $\SL_2(q)$ such that $|\Delta|=7$ and $[K, L]=1$
for every $K,L\in\Delta$  and $K\ne L$. Since $5$ divides $q^2-1$,  Lemma~\ref{2A5} implies that
$\SL_2(q)$ possesses a subgroup isomorphic to $\SL_2(5)$. For every $K\in \Delta$ fix some  $U(K)\le K$
such that $U(K)\cong \SL_2(5)$. Set
 $$
 U=\langle U(K)\mid K\in \Delta \rangle.
 $$
 It follows from the definition that $$U/U_{\mathfrak{S}}\cong \underbrace{\Alt_5\times\dots\times \Alt_5}_{7\text{ times}}$$
 and $U_{\mathfrak{S}}$ is a 2-group. Thus, $U\in\X$. We claim  that $U$ is not conjugate to a subgroup of $H$.
 It is sufficient to show that $|U/U_{\mathfrak{S}}|_5>|H/H_{\mathfrak{S}}|_5$. Indeed, $$|H/H_{\mathfrak{S}}|_5= |\P\Omega_7(2)|_5=5<5^7=|\Alt_5|_5^7=|U/U_{\mathfrak{S}}|_5,$$
a contradiction with $G\in\D_\X$.

 \begin{itemize}
  \item[$(xxx)$]  {\it $G$ is not isomorphic to      $E_8(q)$.}
\end{itemize}

Suppose $G=E_7(q)$ and $H\in\Hall_\X(G)$. By Lemma~\ref{HallSubgroupsOfExceptionalGroups} we have:
\begin{itemize}
  \item  $\{2,3,5,7\}\subseteq \pi\cap \pi(G)\subseteq \pi(q-\varepsilon)$, where $\varepsilon=\pm1$ is such that $4$ divides $q-\varepsilon$, $H$ is a $\pi$-Hall subgroup of a
group
$T\arbitraryext 2\arbitraryext\P\Omega_8^+(2)\arbitraryext 2$, where $T$ is a maximal torus of order $({q-\eta})^8$.
\end{itemize}
This implies that $\P\Omega_8^+(2)\in\X$ and $H/H_{\mathfrak{S}}\cong \P\Omega_8^+(2)\arbitraryext 2$. In $\P\Omega_8^+(2)$ there is a maximal subgroup $\Omega_7(2)$. Therefore, $\Alt_5\in\X$.

Now we argue as in $(xxviii)$. Lemma~\ref{FundamentalSubgroups} implies that  $G$
possesses a collection $\Delta$ of subgroups isomorphic to $\SL_2(q)$ such that $|\Delta|=8$ and $[K, L]=1$
for every $K,L\in\Delta$  and $K\ne L$. Since $5$ divides $q^2-1$,
Lemma~\ref{2A5} implies that  $\SL_2(q)$ possesses a subgroup isomorphic to $\SL_2(5)$.
For every $K\in \Delta$ fix some  $U(K)\le K$ such that $U(K)\cong \SL_2(5)$. Set
 $$
 U=\langle U(K)\mid K\in \Delta \rangle.
 $$
 It follows from the definition that $$U/U_{\mathfrak{S}}\cong \underbrace{\Alt_5\times\dots\times \Alt_5}_{8\text{ times}}$$
 and $U_{\mathfrak{S}}$ is a 2-group. Thus, $U\in\X$. We clam that $U$ is not conjugate to a subgroup of $H$. It is sufficient to show that $|U/U_{\mathfrak{S}}|_5>|H/H_{\mathfrak{S}}|_5$. Indeed, $$|H/H_{\mathfrak{S}}|_5= |\P\Omega_8^+(2)\arbitraryext 2|_5=5^2<5^8=|\Alt_5|_5^8=|U/U_{\mathfrak{S}}|_5,$$  a contradiction with $G\in\D_\X$.

 \begin{itemize}
  \item[$(xxxi)$]  {\it $G$ is not isomorphic to any exceptional group of Lie type.}
\end{itemize}
This statement follows from $(xxvi)$--$(xxx)$.

\subsection{Final proof of the implication $(1)\Rightarrow (2)$}
\begin{itemize}
  \item[$(xxxii)$]  {\it $G$ does not exist.}
\end{itemize}
Indeed, according to the classification of finite simple groups \cite[Theorem~0.1.1]{AschLyoSmSol}, in $(viii)$, $(xv)$, $(xx)$, $(xxv)$, and $(xxxi)$ we have exclude for $G$ all possibilities to be a finite simple group.

Theorem~\ref{DX_means_Dpi} is proven. \qed\medskip
\section{Proofs of Corollaries and Theorem~\ref{Reduktionssatz_main}}

\subsection{
Proofs of Corollaries~\ref{DX_Ext} and~~\ref{ConjMax=ConjSubmax}}

In view of Lemma~\ref{DX_Normal_and_Quot}, in order to prove Corollary~\ref{DX_Ext} it is sufficient to prove that any extension of a $\D_\X$-group by a $\D_\X$-group  is a $\D_\X$-group. Now by Lemma~\ref{DXExtEquiv}, to prove Corollaries~\ref{DX_Ext} and~\ref{ConjMax=ConjSubmax}, it is sufficient to show that if $G\in\D_\X$ is a simple group then $\hat G=\Aut(G)\in\D_\X$.

By Theorem~\ref{DX_means_Dpi} we need to consider two cases: $G\in \X$ and $G\in\D_\X\cap\D_\pi$, where $\pi=\pi(\X)$, and in the last case $G$ is not a $\pi$-group.

In the first case, since $\Aut(G)/\Inn(G)$ is solvable and $\Inn(G)\cong G\in \X$, we conclude that $\hat G$ is $\X$-separable and  $\hat G\in\D_\X$ follows from Lemma~\ref{Bijection}.

In the last case, every $\pi$-Hall subgroup of $G$ is solvable by Lemma~\ref{SolvabilityHallSubgroupsInDpi}. Consequently, the $\X$-subgroups of $G$ are exactly  the solvable $\pi$-subgroups. Since $\Aut(G)/\Inn(G)$ is solvable, the same statement holds for the $\X$-subgroups of $\hat{G}=\Aut(G)$. In particular $\m_\X(\hat G)=\Hall_\pi(\hat G)$. Now Lemma~\ref{DpiExt} implies that $\hat G\in\D_\pi$. Hence the elements of $\m_\X(\hat G)=\Hall_\pi(\hat G)$ are conjugate and $\hat G\in\D_\X$.
\qed\medskip

\subsection{Proof of Corollary~\ref{ClassificationDX}}

Corollary~\ref{DX_Ext} means that $G\in\D_\X$ if and only if every composition factor $S$ of $G$ is a $\D_\X$-group. Now Corollary~\ref{ClassificationDX} immediately follows from Theorem~\ref{DX_means_Dpi} and Lemma~\ref{DpiCrit}. \qed\medskip

\subsection{Proof of Theorem~\ref{Reduktionssatz_main}} Theorem~\ref{Reduktionssatz_main} follows from Corollaries~\ref{DX_Ext}--\ref{ClassificationDX} and \cite[15.4]{Wie4}. \qed\medskip

\subsection{Proof of Corollary~\ref{moduloRadicals}}

We need the following two lemmas which are consequences of Theorem~\ref{Reduktionssatz_main}.

\begin{Lemma}\label{ReductionOfNormalizers} Let
$H\in\m_\X(G)$ for a group $G$. If $A$ is a normal $\D_\X$-subgroup of  $G$ and $\,\overline{\phantom{a}}:G\to G/A$ is the canonical epimorphism  then
$
N_{\overline{G}}(\overline{H})=\overline{N_G(H)}.
$
\end{Lemma}
\noindent{\sc Proof}. We clearly have  $\overline{N_G(H)}\le N_{\overline{G}}(\overline{H})$.
Suppose $x\in N_G(HU)$ (equivalently, ${\overline{x}\in N_{\overline{G}}(\overline{H})}$). Then $H^xA=HA$. The $\X$-Reduktionssatz for $A$
implies that there is $a\in A$ such that $H^x=H^a$, i.\,e. $xa^{-1}\in N_G(H)$,  $x\in N_G(H)A$, and $\overline{x}\in\overline{N_G(H)}$.
\qed\medskip

\begin{Lemma}\label{three}
  Let $G\,{\trianglelefteqslant\trianglelefteqslant}\, G^*$ and let $H=K\cap G$ for some  $K\in\m_\X(G^*)$.
  Suppose that $A\in\D_\X$ is a normal subgroup of $G^*$ and $\,\overline{\phantom{a}}:G^*\to G^*/A$ is the canonical epimorphism.
  Then $\overline{H}=\overline{K}\cap\overline{G}$.
\end{Lemma}
\noindent{\sc Proof}. Clearly, $\overline{H}=\overline{K\cap G}\le \overline{K}\cap \overline{G}$. Suppose that
$\overline{H}< \overline{K}\cap \overline{G}$. Since $G\,{\trianglelefteqslant\trianglelefteqslant}\, G^*$,
we have $H\,{\trianglelefteqslant\trianglelefteqslant}\, K$, $\overline{H} \,{\trianglelefteqslant\trianglelefteqslant}\, \overline{K}$,
and $\overline{H} \,{\trianglelefteqslant\trianglelefteqslant}\, \overline{K}\cap \overline{G}$.

Since $\overline{H}<\overline{K}\cap \overline{G}$, the index of $\overline{H}$ in $N_{\overline{G}}(\overline{H})$ is divisible by a prime
$p\in \pi(\X)$. By Lemma~\ref{ReductionOfNormalizers} we have
$
N_{\overline{G}}(\overline{H})=\overline{N_G(H)}.
$
Now
$$
p\mid \left|N_{\overline{G}}(\overline{H}):\overline{H}\right|=\left|\frac{N_G(H)A}{HA}\right|=\frac{|N_G(H)|}{|N_A(H)|}:\frac{|H|}{|H\cap A|}=\frac{|N_G(H):H|}{|N_A(H):(H\cap A)|}.
$$
Therefore, $p\mid \left|N_G(H)/H\right|$ contrary to Lemma \ref{WH-strong}.
\qed\medskip

Now we prove Corollary~\ref{moduloRadicals}. Denote by $\phi$ the canonical epimorphism $G\rightarrow G/N$.
We need to show that  $H^\phi\in\sm_\X(G^\phi)$. We can assume that there exists a group $G^*$ such that $G\,{\trianglelefteqslant\trianglelefteqslant}\, G^*$ and $H=G\cap K$ for some $K\in\m_\X(G^*)$. Denote by $A$ the normal closure $\langle N^x\mid {x\in G^*}\rangle$ of $N$ in~$G^*$. As $N^x\trianglelefteqslant G^x\, {\trianglelefteqslant\trianglelefteqslant}\, G^*$ for every $x\in G^*$ and $\mathfrak{F}$ is a Fitting class\footnote{For a Fitting class $\mathfrak{F}$, if $N_1,\dots N_s$ are subnormal $\mathfrak{F}$-subgroups of $G$ then $\langle N_1,\dots N_s\rangle\in \mathfrak{F}$ \cite[II(2.8)]{DH}.}, we conclude that $A\in\mathfrak{F}$. Moreover, $G\cap A\,{\trianglelefteqslant\trianglelefteqslant}\, A$. Consequently, $G\cap A\in\mathfrak{F}$. Since $G\cap A\trianglelefteqslant G$, we have
$$N\leq G\cap A\leq G_{\mathfrak{F}}=N.$$
Therefore, $G\cap A=N$.
Consider the restriction $\tau:G\rightarrow G^*/A$ to $G$ of the canonical epimorphism $G^*\rightarrow G^*/A$. Then the kernel $G\cap A$ of $\tau$ coincides with $N$. By the homomorphism theorem, there exists an injective homomorphism
$\psi: G^\phi=G/A \rightarrow G^*/A$ such that the following diagram is commutative:
$$\xymatrix{
  G \ar[d]_{\phi} \ar[r]^{\tau} &     G^*/A   \\
 G^\phi \ar[ur]_{\psi}                     }$$
Then  $G^{\phi\psi}=G^\tau={GA/A}\,{\trianglelefteqslant\trianglelefteqslant}\, {G^*/A}.$
Moreover, by Lemma~\ref{three} we have $$H^{\phi\psi}=H^\tau={HA/A}=(G\cap K)A/A=(GA/A)\cap (KA/A)=G^{\phi\psi}\cap (KA/A),$$  where $KA/A\in\m_\X(G^*/A)$  by Theorem~\ref{Reduktionssatz_main}. Thus $H^\phi\in\sm_\X(G^\phi)$ by the definition of an $\X$-sub\-maxi\-mal subgroup.

In order to prove the inequality $\overline{k}\geqslant k$, it is sufficient to prove that if $H,H_1\in\sm_\X(G)$ and $HN/N$ and $H_1N/N$ are conjugate in $G/N$, then $H$ and $H_1$ are conjugate in~$G$. Without loss of generality, we can assume that $HN/N=H_1N/N$ and $HN=H_1N$. By Corollary~\ref{DX_Ext} we have $G_0:=HN\in\D_\X$. Moreover,   $$H\cap N\in \sm_\X(N)=\m_\X(N)=\Hall_\X(N)$$ by Corollary~\ref{ConjMax=ConjSubmax} and Lemma~\ref{HallXSubgroups}. Therefore, $H\in\Hall_\X(G_0)=\m_\X(G_0)$. Similarly,  ${H_1\in\m_\X(G_0)}$. Hence, $H$ and $H_1$ are conjugate in~$G_0$.
 \qed\medskip

\section*{Appendix}

In Conditions I--VII below, $S$ is a finite simple group, $\X$ is a complete class of finite groups, and $\pi=\pi(\X)$.

\medskip

\noindent{\bf Condition I.}  We say that $(S,\X)$ satisfies Condition I
if $S\in\X$ or $|\pi\cap\pi(S)|\le 1$.

\medskip

\noindent{\bf Condition II.} We say that $(S,\X)$  satisfies
Condition  II if one of the following cases holds.
\begin{itemize}
\item[$(1)$]  $S\simeq M_{11}$ and $\pi\cap\pi(S)=\{5,11\}$;
\item[$(2)$]  $S\simeq M_{12}$  and $\pi\cap\pi(S)=\{5,11\}$;
\item[$(3)$]  $S\simeq M_{22}$ and $\pi\cap\pi(S)=\{5,11\}$;
\item[$(4)$]  $S\simeq M_{23}$  and $\pi\cap\pi(S)$
coincide with one of the following sets $\{5,11\}$ and $\{11,23\}$;
\item[$(5)$]  $S\simeq M_{24}$  and $\pi\cap\pi(S)$
coincide with one of the following sets $\{5,11\}$ and $\{11,23\}$;
\item[$(6)$]  $S\simeq J_1$ and $\pi\cap\pi(S)$ coincide with one of the following sets $\{3,5\}$,
$\{3,7\}$,  $\{3,19\}$, and  $\{5,11\}$;
\item[$(7)$]  $S\simeq J_4$ and $\pi\cap\pi(S)$ coincide with one of the following sets  $\{5,7\}$,
$\{5,11\}$, $\{5,31\}$, $\{7,29\}$, and $\{7,43\}$;
\item[$(8)$]  $S\simeq O'N$ and $\pi\cap\pi(S)$
coincide with one of the following sets $\{5,11\}$ and $\{5,31\}$;
\item[$(9)$]  $S\simeq Ly$ and $\pi\cap\pi(S)=\{11,67\}$;
\item[$(10)$]  $S\simeq Ru$ and $\pi\cap\pi(S)=\{7,29\}$;
\item[$(11)$]  $S\simeq Co_1$ and $\pi\cap\pi(S)=\{11,23\}$;
\item[$(12)$]  $S\simeq Co_2$  and $\pi\cap\pi(S)=\{11,23\}$;
\item[$(13)$]  $S\simeq Co_3$ and $\pi\cap\pi(S)=\{11,23\}$;
\item[$(14)$]  $S\simeq M(23)$  and $\pi\cap\pi(S)=\{11,23\}$;
\item[$(15)$]  $S\simeq M(24)'$  and $\pi\cap\pi(S)=\{11,23\}$;
\item[$(16)$]  $S\simeq B$  and $\pi\cap\pi(S)$ coincide with one of the following sets $\{11,23\}$ and  $\{23,47\}$;
\item[$(17)$]  $S\simeq M$ and $\pi\cap\pi(S)$ coincide with one of the following sets $\{23,47\}$ and  $\{29,59\}$.
\end{itemize}

\medskip

\noindent{\bf Condition III.} Let $S$ be isomorphic to a group of Lie type over the field $\F_q$ of characteristic $p\in\pi$ and let
$\tau=(\pi\cap\pi(S))\setminus\{p\}$. We say that $(S,\X)$
satisfies Condition~III if $\tau\subseteq\pi(q-1)$ and every
prime in $\pi$ does not divide the order of the Weyl group of~$S$.

\medskip

In order to formulate Conditions IV and V, we need the following notation. If $r$~is an odd prime and  $q$~is an integer not divisible by $r$, then
$e(q,r)$ is the smallest
 positive integer  $e$ with $q^e \equiv 1 \pmod r$.

\medskip
\noindent{\bf Condition IV.} Let $S$ be  isomorphic to   a group of Lie type
with the base field $\F_q$ of  characteristic $p$.
Let
$2,p\not\in\pi$. Denote by $r$ the minimum in $\pi\cap\pi(S)$ and put ${\tau=(\pi\cap\pi(S))\setminus\{r\}}$ and  $a=e(q,r)$. We say that
$(S,\X)$  satisfies Condition IV if there exists $t\in\tau$
with $b=e(q,t)\ne a$ and one of the following statements holds.

\begin{itemize}

\item[$(1)$]\noindent  $S\simeq A_{n-1}(q)$, $a=r-1$, $b=r$, $(q^{r-1}-1)_r=r$,
$\left[\frac{n}{r-1}\right]=\left[\frac{n}{r}\right]$,
and $e(q,s)=b$ for every $s\in\tau$;

\item[$(2)$]  $S\simeq A_{n-1}(q)$, $a=r-1$, $b=r$, $(q^{r-1}-1)_r=r$,
$\left[\frac{n}{r-1}\right]=\left[\frac{n}{r}\right]+1$,
$n\equiv -1 \pmod r$, and $e(q,s)=b$ for every $s\in\tau$;

\item[$(3)$]  $S\simeq {^2A}_{n-1}(q)$, $r\equiv 1 \pmod 4$, $a=r-1$, $b=2r$, $(q^{r-1}-1)_r=r$,
$\left[\frac{n}{r-1}\right]=\left[\frac{n}{r}\right]$
and $e(q,s)=b$ for every $s\in\tau$;

\item[$(4)$]  $S\simeq {^2A}_{n-1}(q)$, $r\equiv 3 \pmod 4$, $a=\frac{r-1}{2}$,
$b=2r$, $(q^{r-1}-1)_r=r$,
$\left[\frac{n}{r-1}\right]=\left[\frac{n}{r}\right]$
and  $e(q,s)=b$ for every $s\in\tau$;

\item[$(5)$]  $S\simeq {^2A}_{n-1}(q)$, $r\equiv 1 \pmod 4$, $a=r-1$, $b=2r$, $(q^{r-1}-1)_r=r$,
$\left[\frac{n}{r-1}\right]=\left[\frac{n}{r}\right]+1$,
$n\equiv -1 \pmod r$ and $e(q,s)=b$ for every $s\in\tau$;

\item[$(6)$]  $S\simeq {^2A}_{n-1}(q)$, $r\equiv 3 \pmod 4$, $a=\frac{r-1}{2}$,
$b=2r$, $(q^{r-1}-1)_r=r$,
$\left[\frac{n}{r-1}\right]=\left[\frac{n}{r}\right]+1$,
$n\equiv -1 \pmod r$ and
$e(q,s)=b$ for every  $s\in\tau$;

\item[$(7)$]  $S\simeq {^2D}_n(q)$, $a\equiv 1 \pmod 2$, $n=b=2a$ and for every $s\in\tau$ either $e(q,s)=a$ or
$e(q,s)=b$;

\item[$(8)$]  $S\simeq {^2D}_n(q)$, $b\equiv 1 \pmod 2$, $n=a=2b$ and for every $s\in\tau$ either $e(q,s)=a$ or
$e(q,s)=b$.
\end{itemize}

\medskip

\noindent{\bf Condition V.} Let $S$ be  isomorphic to  a group of Lie type
with the base field $\F_q$ of characteristic $p$. 
Suppose,
$2,p\not\in\pi$.  Let $r$ be the minimum in $\pi\cap\pi(S)$, let ${\tau=(\pi\cap\pi(S))\setminus\{r\}}$, and let  $c=e(q,r)$.  We say that
$(S,\X)$ satisfies Condition V if $e(q,t)=c$ for every
$t\in\tau$ and one of the following statements holds.
\begin{itemize}
\item[$(1)$]  $S\simeq A_{n-1}(q)$ and $n<cs$ for every $s\in\tau$;

\item[$(2)$]  $S\simeq{^2A}_{n-1}(q)$, $c\equiv 0 \pmod 4$ and $n<cs$ for every $s\in\tau$;

\item[$(3)$]  $S\simeq{^2A}_{n-1}(q)$, $c\equiv 2 \pmod 4$ and ${2n<cs}$ for every $s\in\tau$;

\item[$(4)$]  $S\simeq{^2A}_{n-1}(q)$, $c\equiv 1 \pmod 2$ and ${n<2cs}$ for every $s\in\tau$;

\item[$(5)$]  $S$ is isomorphic to one of the groups $B_n(q)$, $C_n(q)$, or ${^2D}_n(q)$,  $c$
is odd and $2n<cs$ for every $s\in\tau$;

\item[$(6)$]   $S$ is isomorphic to one of the groups $B_n(q)$, $C_n(q)$, or $D_n(q)$, $c$
is even and $n<cs$ for every $s\in\tau$;

\item[$(7)$]  $S\simeq D_n(q)$,  $c$ is even and  $2n\leq cs$ for every $s\in\tau$;

\item[$(8)$]  $S\simeq {^2D}_n(q)$,  $c$  is odd and $n\leq cs$ for every $s\in\tau$;

\item[$(9)$]  $S\simeq {^3D}_4(q)$;

\item[$(10)$]  $S\simeq E_6(q)$,  and if $r=3$ and $c=1$ then $5,13\not\in\tau$;

\item[$(11)$]  $S\simeq {^2E}_6(q)$, and if $r=3$ and  $c=2$ then $5,13\not\in\tau$;

\item[$(12)$]  $S\simeq E_7(q)$, if $r=3$ and $c\in\{1,2\}$ then $5,7,13\not\in\tau$,
and if $r=5$ and $c\in\{1,2\}$ then
$7\not\in\tau$;

\item[$(13)$]  $S\simeq E_8(q)$, if
 $r=3$ and $c\in\{1,2\}$ then $5,7,13\not\in\tau$, and if $r=5$ and
$c\in\{1,2\}$ then $7,31\not\in\tau$;

\item[$(14)$]  $S\simeq G_2(q)$;

\item[$(15)$]  $S\simeq F_4(q)$, and if $r=3$ and $c=1$ then $13\not\in\tau$.
\end{itemize}

\medskip

\noindent{\bf Condition VI.} We say that $(S,\X)$ satisfies Condition VI
if one of the following statements holds.

\begin{itemize}
\item[$(1)$]  $S$ is isomorphic to ${^2B}_2(2^{2m+1})$ and  $\pi\cap\pi(S)$ is contained in one of the sets
$$\pi(2^{2m+1}-1), \,\,\, \pi(2^{2m+1}\pm 2^{m+1}+1);$$

\item[$(2)$]  $S$ is isomorphic to ${^2G}_2(3^{2m+1})$ and  $\pi\cap\pi(S)$ is contained in one of the sets
$$\pi(3^{2m+1}-1)\setminus\{2\}, \,\,\, \pi(3^{2m+1}\pm 3^{m+1}+1)\setminus\{2\};$$

\item[$(3)$]  $S$ is isomorphic to ${^2F}_4(2^{2m+1})$ and $\pi\cap\pi(S)$ is contained in one of the sets
$$\pi(2^{2(2m+1)}\pm 1), \,\,\, \pi(2^{2m+1}\pm 2^{m+1}+1),$$
$$\pi(2^{2(2m+1)}\pm 2^{3m+2}\mp2^{m+1}-1),\,\,\, \pi(2^{2(2m+1)}\pm
2^{3m+2}+2^{2m+1}\pm2^{m+1}+1).$$
\end{itemize}

\medskip

\noindent{\bf Condition VII.} Let $S$ be  isomorphic to  a group of Lie type with the base field
$\F_q$ of characteristic $p$. Suppose that  $2\in \pi$ and  $3, p\not\in\pi$. Put
 $\tau=(\pi\cap\pi(S))\setminus\{
2\}$ and  $\varphi=\{t\in\tau\mid t \text{ is a Fermat number}\}$. We say
that $(S,\X)$  satisfies Condition VII if
$\tau\subseteq\pi(q-\varepsilon)$, where the number $\varepsilon=\pm 1$ is  such that $4$ divides $q-\varepsilon$, and one of the following statements
holds.
\begin{itemize}
\item[$(1)$] $S$ is isomorphic to either ${A}_{n-1}(q)$ or  ${}^2{A}_{n-1}(q)$, $s>n$
 for every $s\in\tau$, and $t>n+1$ for every $t\in\varphi$;
\item[$(2)$] $S\simeq {B}_n(q)$, and $s>2n+1$
for every $s\in\tau$;
\item[$(3)$] $S\simeq {C}_n(q)$, $s>n$ for every
  $s\in\tau$, and $t>2n+1$ for every  $t\in\varphi$;
\item[$(4)$] $S$ is isomorphic to either ${D}_n(q)$ or  ${}^2{D}_n(q)$,  and $s>2n$
for every $s\in\tau$;
\item[$(5)$] $S$ is isomorphic to either ${G}_2(q)$ or ${^2G}_2(q)$, and $7\not\in\tau$;
\item[$(6)$] $S\simeq {F}_4(q)$ and $5,7\not\in\tau$;
\item[$(7)$] $S$ is isomorphic to either $E_6(q)$ or ${}^2E_6(q)$, and $5,7\not\in\tau$;
\item[$(8)$] $S\simeq {E}_7(q)$ and $5,7,11\not\in\tau$;
\item[$(9)$] $S\simeq {E}_8(q)$ and $5,7,11,13\not\in\tau$;
\item[$(10)$] $S\simeq {^3D}_4(q)$ and $7\not\in\tau$.
\end{itemize}

%
%
%
%
%
%
%
%
%
%
%
%
%
%
%
%
%
%
%
%
%
%
%
%
%
%
%
%
%
%
%
%


\smallskip

\vspace{3cm}


\small
\begin{thebibliography}{100}


\bibitem{Asch1} M.\,Aschbacher, Characterization of Chevalley groups over fields of odd order, Ann. Math., 106 (1977), 353--468.

\bibitem{Asch2} M.\,Aschbacher, On finite groups of Lie type and odd characteristic, J. Algebra, 66:2 (1980), 400--424.

\bibitem{AschLyoSmSol} M.\,Aschbacher, R.\,Lyons, S.\,D.\,Smith, R.\,Solomon, The classification of finite simple groups. Groups of characteristic 2 type. Mathematical Surveys and Monographs, 172. American Mathematical Society, Providence, RI, 2011. xii+347 pp.


\bibitem{Bray}
{J.\,N.\,Bray, D.\,F.\,Holt, C.\,M.\,Roney-Dougal}, The Maximal Subgroups of the Low-Dimensional Finite Classical Groups. Cambridge: Cambridge Univ. Press, 2013. 438~p.

\bibitem{CF} R.\,Carter, P.\,Fong, The Sylow $2$-subgroups of the finite classical groups, J. Algebra, 1:1 (1964), 139--151.


\bibitem{CLSS}
{A.\,M.\,Cohen, M.\,W.\,Liebeck, J.\,Saxl,  G.\,M.\,Seitz},
The local maximal subgroups of exceptional groups of Lie type, finite and algebraic,
{Proc. London Math. Soc.} Ser. III, {64} (1992), N1, 21--48.




\bibitem{Atlas}
{J.\,H.\,Conway, R.\,T.\,Curtis, S.\,P.\,Norton, R.\,A.\,Parker, R.\,A.\,Wilson,}
Atlas of Finite Groups. Oxford: Clarendon Press, 1985. 252~p.




\bibitem{ChuShem} S.\,A.\,Chunikhin, L.\,A.\,Shemetkov, Finite groups, J. Soviet Math. 1:3 (1973),
291--332.


\bibitem{DH}
{K.\,Doerk, T.\,Hawks,}
{Finite Soluble Groups},
Berlin, New York, Walter de Gruyter, 1992.




\bibitem{Praeger} {D.\,Easdown, C.E.\,Praeger}, On minimal faithful permutation representations of finite groups, Bulletin Australian Mathematical Society, \textbf{38}(1988), 207--220.



\bibitem{FeitThompsonOddOrder}
{ W.\,Feit, J.\,G.\.Thompson},
Solvability of groups of odd order,
{ Pacif. J. Math.}, {\bf 13} (1963), N~3, 775--1029.

\bibitem{Galois}  \'{E}.\,Galois, M\'{e}moire sur les conditions de r\'{e}solubilit\'{e} des \'{e}quations par radicaux, J. Math. Pures Appl. (Liouville) 11 (1846), 417--433.


\bibitem{Glauberman} {G.\,Glauberman}, Factorization in local subgroups of finite groups.  Conf. Ser. Math., 33, Amer. Math. Soc., Providence, RI, 1976.


\bibitem{Gross1} F.\,Gross, On a conjecture of Philip Hall, Proc. Lond. Math. Soc. (3), 52:3 (1986), 464--494.

\bibitem{Gross3} F.\,Gross, Conjugacy of odd order Hall subgroups, Bull. Lond. Math. Soc., 19:4 (1987), 311--319.

\bibitem{Gross4} F.\,Gross, Hall subgroups of order not divisible by 3, Rocky Mount. J. Math., 23:2 (1993), 569--591.

\bibitem{Gross2} F.\,Gross, Odd order Hall subgroups of the classical linear groups, Math. Z., 220:3 (1995), 317--336.



\bibitem{GuoBook}
{W.\,Guo,}
{The Theory of Classes of Groups},
Beijing, New York, Kluwer Acad. Publ., 2006.

\bibitem{GuoBook1}
{W.\,Guo,}
{Structure Theory of Canonical Classes of Finite Groups}, Berlin, Springer, 2015.



\bibitem{GR1} W.\,Guo, D.\,O.\,Revin,   On maximal and submaximal ${\mathfrak X}$-subgroups, Algebra and logic, 56:1 (2018),  9--28.

\bibitem{GR2} W.\,Guo, D.\,O.\,Revin,   Conjugacy of maximal and submaximal $\mathfrak{X}$-subgroups, Algebra and Logic, 57:3 (2018), 169--181.

\bibitem{GR_SubmaxMinNonSolv} W.\,Guo, D.\,O.\,Revin,   Classification and properties of the $\pi$-submaximal
subgroups in minimal nonsolvable groups,  Bulletin of Mathematical Sciences, 8:2 (2018), 325--351.

\bibitem{GR_surv} W.\,Guo, D.\,O.\,Revin,  Pronormality and submaximal $\mathfrak{X}$-subgroups in finite groups, Communications in Mathematics and Statistics, 6:3 (2018), 289--317.

\bibitem{GRV} W.\,Guo, D.\,O.\,Revin, E.\,P.\,Vdovin,  Confirmation for Wielandt's conjecture, J. Algebra, 434 (2015), 193--206.


\bibitem{Hall} P. Hall,  Theorems like Sylow`s,
 Proc. London Math. Soc., 6:22 (1956), 286--304.



\bibitem{Hart}
{ B.\,Hartley, }
A theorem of Sylow type for a finite groups,
{ Math.\ Z.}, { 122:4} (1971),  223--226.




  \bibitem{Holt} {D.\,F.\,Holt}, Representing quotients of permutation groups, Quarterly Journal of Mathematics, {48}, No. 2 (1997), 347--350.



\bibitem{Jordan0}  C.\,M.\,Jordan,  Commentaire sur le M\'{e}moire de Galois. Comptes rendus 60 (1865), 770--774.
\bibitem{Jordan} C.\,M.\,Jordan, Trait\'{e} des substitutions et des \'{e}quations alg\'{e}briques, Paris: Gauthier-Villars, 1870.


 \bibitem{KL}
{P.\,B.\,Kleidman,  M.\,Liebeck,} The subgroup structure of the finite classical groups. Cambridge: Cambridge Univ. Press, 1990. 303~p.


\bibitem{Manz}
  N.\,Ch.\,Manzaeva, Heritability of the property $\mathcal D_\pi$ by overgroups of $\pi$-Hall subgroups in the case where $2\in\pi$, Algebra and Logic, 53:1 (2014), 17--28.




\bibitem{Kour} V.\,D.\,Mazurov and E.\,I.\,Khukhro (eds.), The Kourovka notebook. Unsolved
problems in group theory, 17th ed., Institute of Mathematics, Siberian Branch of
RAS, Novosibirsk 2010.

\bibitem{MR} V.\,D.\,Mazurov, D.\,O.\,Revin, On the Hall $D_\pi$-property for finite groups, Siberian Math. J., 38:1 (1997), 106--113.


\bibitem{HallChev} D. O. Revin, Hall $\pi$-subgroups of finite Chevalley groups whose characteristic belongs to~$\pi$, Siberian Advances in Mathematics, 1999, 9:2, 25--71.

 \bibitem{DpiCl} D. O. Revin, The $D_\pi$-property in a class of finite Groups, Algebra and Logic, 41:3 (2002), 187--206.



\bibitem{R2}  D.\,O. Revin, The $D_\pi$ property of finite groups in the case $2\notin\pi$, Proceedings of the Steklov Institute
 of Mathematics, 257:suppl.1 (2007), S164--S180.



\bibitem{R4}  D.\,O. Revin, The $D_\pi$-property in finite simple groups, Algebra and Logic, 47:3 (2008), 210--227.



\bibitem{R5} D.\,O. Revin, The $D_\pi$-property of linear and unitary groups, Siberian Math. J., 49:2 (2008), 353--361.


\bibitem{Epimax} D. O. Revin, submaximal and epimaximal $\mathfrak{X}$-subgroups, Algebra and Logic, 58:6 (2019), 475--479.


\bibitem{RSV} D. O. Revin, S. V. Skresanov, A. V. Vasilev, The Wielandt--Hartley theorem for submaximal $\mathfrak{X}$-subgroups, Monatshefte f\"{u}r Mathematik, 193:1 (2020),  143--155.


 \bibitem{Hall3'} D.\,O.\,Revin, E. P. Vdovin, Hall subgroups of finite groups, Contemporary Mathematics, 402 (2006), 229--265.




 \bibitem{NumbCl} D.\,O.\,Revin, E.\,P.\,Vdovin,  On the number of classes of conjugate Hall subgroups in finite simple groups, J. Algebra, 324:12 (2010),  3614--3652.


\bibitem{RZ1} D.\,O.\,Revin, A.\,V.\,Zavarnitsine, On the behavior of $\pi$-submaximal subgroups under homomorphisms, Communications in Algebra, 48:2 (2020), 702--707.

\bibitem{RZ2} D.\,O.\,Revin, A.\,V.\,Zavarnitsine, The behavior of $\pi$-submaximal subgroups under homomorphisms with  $\pi$-separable kernels, Siberian Electronic Mathematical Reports, 17 (2020), accepted, see 	\texttt{arXiv:2006.09752}.


\bibitem{ShemBook}  L.\,A.\, Shemetkov, Formations of Finite Groups (in Russian), Nauka, Moscow (1978).


\bibitem{Shem7} L.\,A.\,Shemetkov, Two
directions in the development of the theory of non-simple finite groups, Russian
Math. Surveys 30:2 (1975), 185--206.


\bibitem{Suz} M.\,Suzuki, { Group Theory} II,  Springer-Verlag, New
York--Berlin--Heidelberg--Tokyo, 1986.



\bibitem{Thomp}
{J.\,G.\,Thompson},
 Nonsolvable finite groups all of whose local subgroups are solvable,  Bull. Amer. Math. Soc., 74 (1968), 383--437.


\bibitem{VMR}
  E.\,P.\,Vdovin, N.\,Ch.\,Manzaeva, D.\,O.\,Revin, On the heritability of the property $D_\pi$ by subgroup, Proceedings of the Steklov Institute of Mathematics, 279:Suppl.~1 (2012), 130--138.

\bibitem{MRV}
 E.\,P.\,Vdovin, N.\,Ch.\,Manzaeva, D.\,O.\,Revin, On the heritability the Sylow $\pi$-theorem by subgroups, Sb. Math., 211:3 (2020),
309--335.


\bibitem{Surv} E.\,P.\,Vdovin, D.\, O.\,Revin, Theorems of Sylow type,
    Russian Math. Surveys, 66:5 (2011), 829--870.




\bibitem{Weir1955} A.J. Weir, Sylow $p$-subgroups of the
classical groups over finite fields with characteristic prime to $p$,
Proceedings of the American Mathematical Society 6:4 (1955), 529--533.



\bibitem{Wie}
{ H. Wielandt},
Zum Satz von Sylow,
{ Math.\ Z.},  { 60} (1954), N4. 407--408.


\bibitem{Wie1}
{ H.\,Wielandt},
Entwicklungslinien in der Strukturtheorie der endlichen Gruppen,
{ Proc. Intern. Congress Math., Edinburg, 1958.} London: Cambridge
Univ. Press, 1960, 268--278.


\bibitem{Wie5} H.\,Wielandt, Arithmetische Struktur und Normalstruktur endlicher Gruppen,
Conv. Internaz. di Teoria dei Gruppi Finiti e Applicazioni (Firenze 1960), Edizioni
Cremonese, Roma 1960,  56--65.


\bibitem{Wie4}
{H.\,Wielandt}
Zusammengesetzte Gruppen endlicher Ordnung,
{ Vorlesung an der Universit\"at T\"ubingen im Wintersemester 1963/64}.
Helmut Wielandt: Mathematical Works, Vol. 1,
Group theory (ed. B. Huppert and H. Schneider, de Gruyter, Berlin, 1994), 607--655.


\bibitem{WieCanb}
{ H.\,Wielandt},
On the structure of composite groups, Proc. Internat. Conf. Theory of Groups, Austral. Nat. Univ. Canberra, August 1965,
pp. 379-388. Gordon and Breach Science Publishers, Inc., New York 1967.




\bibitem{Wie3}
{H.\,Wielandt},
Zusammengesetzte Gruppen: H\"older Programm heute,
{ The Santa Cruz conf. on finite groups, Santa Cruz, 1979}. Proc.
Sympos. Pure
Math., { 37}, Providence RI: Amer. Math. Soc., 1980, 161--173.


\bibitem{Wie_diary}
{H.\,Wielandt}, Tageb\"{u}cher, D17, 
{{(Mathematical Diary XVII)}}, 1980, available in photocopied and transcribed form at \texttt{https://www3.math.tu-berlin.de/numerik/Wielandt/index\_en.html}


\bigskip
\end{thebibliography}
\end{document}